\documentclass[12pt,a4paper]{article}
\usepackage{amsmath}
\usepackage{amssymb}
\advance\topmargin by -2.2cm                    
\newtheorem{theo}{\indent Theorem}[section]
\newtheorem{prop}[theo]{\indent Proposition}

\newtheorem{rem}[theo]{\indent Remark}

\newtheorem{defin}[theo]{\indent Definition}
\newtheorem{cor}[theo]{\indent Corollary}

\newtheorem{ex}[theo]{\indent Example}
\newtheorem{ass}[theo]{\indent Assumption}

\newenvironment{proof}{\noindent {\bf Proof }}
{\hfill $\bullet$ \vspace{0.25cm}}

\def\RR{\mathbb{R}}
\def\1{\mathbf{1}}
\def \E{\mathbb{E}}
\def \P{\mathbb{P}}

\def\ZZ{\mathbb{Z}}

\def \ge{\varepsilon}
\def \FF{I\!\!F}
\def \GG{C\!\!\!\!G}
\begin{document}
\newlength{\breite}
\settowidth{\breite}{Web: http://www.mathematik.uni-mainz.de/Stochastik/Loecherbache}

\title{Deviation inequalities for centered additive functionals of recurrent Harris processes having general state space}

\author{Eva 
{\sc L\"ocherbach}\footnote{LAMA UMR-CNRS 8050,
Facult\'e de Sciences et Technologie,
Universit\'e Paris-Est Val de Marne, 61 avenue du G\'en\'eral de Gaulle, 94010 Cr\'eteil
Cedex, France. E-mail: {\tt
locherbach@univ-paris12.fr}} , Dasha {\sc Loukianova}\footnote{D\'epartement de Math\'ematiques, Universit\'e d'Evry-Val d'Essonne, Bd Fran\c{c}ois Mitterrand, 91025 Evry Cedex, France. E-mail: {\tt 
 dasha.loukianova@univ-evry.fr
} }
}

\maketitle
\def\abstractname{Abstract}
\begin{abstract}
Let $X$ be a Harris recurrent strong Markov process in continuous time
 with general Polish state space $E,$ having invariant measure $\mu .$  In this
paper we use the regeneration method to derive non asymptotic deviation bounds for
$$P_{x}\left (|\int_0^tf(X_s)ds|\geq t^{\frac12 + \eta} \ge\right)$$
in the positive recurrent case, for nice functions $f$ with $\mu (f)
=0 $ ($f$ must be a charge). We generalize these bounds to the fully null-recurrent case in the moderate deviations regime. We
obtain a Gaussian contentration bound for all functions $f$ which are a charge. The rate of convergence is expressed in
terms of the deterministic equivalent of the process. The main
ingredient of the proof is Nummelin splitting in continuous time
which allows to introduce regeneration times for the process on an enlarged state space.

\end{abstract}

{\it Key words} : Harris recurrence, Nummelin splitting, continuous
time Markov processes, special functions, additive functionals, 
deviation inequalities, deterministic equivalent
for additive functionals.
\\

{\it MSC 2000}  : 60 J 55, 60 J 35, 60 F 10, 62 M 05

\section{Introduction}
Consider a Harris recurrent strong Markov process $X = (X_t)_{t\geq
0 }$ in continuous time with invariant measure $\mu ,$ taking values in a Polish space
$E.$

If the total mass of $\mu$ is finite, $X$ is called positive recurrent, null-recurrent otherwise.
In the case of positive recurrence it is well known that for certain functions $f$
 such that $\mu(f) = 0 $ we have a central limit theorem for 
\begin{equation}\label{eq:clt}
\left(  \frac{1}{\sqrt{n}} \int_0^{nt} f(X_s) ds \right)_{t \geq 0} 
\end{equation}
as $n$ goes to $\infty,$ see for instance H\"opfner and L\"ocherbach \cite{ams} and Touati \cite{touati90}.

In the null-recurrent case, the re-normalization $\sqrt{n} $ of (\ref{eq:clt}) has to be changed. For that sake we have to consider what is called {\it deterministic
equivalent} of additive functionals. The deterministic equivalent has been introduced for Markov chains
by Chen \cite{chen} and then been generalized to the context of 
continuous time diffusion models by Loukianova and Loukianov \cite{dashaoleg}
and to any continuous time recurrent Markov process by L\"ocherbach and 
Loukianova in \cite{dashaeva}. 
It is a deterministic function $ t \mapsto v (t) $ such that $v (t) \to \infty$ as $t \to \infty $ and such that for any integrable additive functional $A_t ,$
$$ \lim_{M \to \infty } \lim \inf_{t \to \infty } P_\pi ( 1/M \le A_t / v(t) \le M ) = 1 $$
for any initial measure $\pi .$ $v(t)$ can be defined as follows. Take any fixed positive special function $g$ of the process having $\mu (g) > 0$  (see definition \ref{specialdef} below for the exact definition of special functions, for strong Feller processes, any bounded function having compact support is special) and define 
$$ v(t) := E_\eta ( \int_0^t g(X_s) ds ) ,$$
where $\eta $ is an arbitrary initial measure. Then the strong Chacon-Ornstein theorem implies that for any other special function $g'$ and any other initial measure $\eta' , $ 
$$ \lim_{t \to \infty } \frac{E_\eta ( \int_0^t g(X_s) ds )}{E_{\eta '} ( \int_0^t g' (X_s) ds )}= \frac{\mu (g)}{\mu (g')} .$$
Hence the deterministic equivalent is unique up to a constant in the sense that for two choices of the deterministic equivalent, $v$ and $v'$, we have that $ \lim_{t \to \infty }  v(t)/v'(t) = c ,$ where $c$ is a positive constant. In {\it regular} models, $v(t) \sim t^{\alpha } l(t),$ where $l$ is a function that varies slowly at infinity. For example, for Brownian motion in dimension one, we have $\alpha = 1/2 .$ 
The explosion rate $v_t$ is in general slower than the
ergodic (positive recurrent) rate $t.$ 

In the null-recurrent regular case, we have convergence in law of
\begin{equation}\label{eq:clt2}
\left(  \frac{1}{\sqrt{v(n)}} \int_0^{nt} f(X_s) ds \right)_{t \geq 0} 
\end{equation} 
to $ B \circ W^{\alpha } ,$ where $B$ is a one-dimensional Brownian motion and $W^\alpha $ the Mittag-Leffler process of index $\alpha ,$ $ 0 < \alpha \le 1 ,$  for certain functions $f$ having $\mu (f)= 0,$ see Touati (\cite{touati90}) and H\"opfner and L\"ocherbach (\cite{ams}) for similar results for martingale additive
functionals. In the fully null-recurrent case (i.e. null-recurrent, but not
regular) no convergence in law holds true, which is a consequence of the famous theorem of Darling-Kac.   

Continuing in the spirit of the limit theorems of (\ref{eq:clt}) and (\ref{eq:clt2}), the aim of this paper is to study  
deviation inequalities for additive functionals $ \int_0^t f(X_s) ds $  for a certain class of
centered functions $f$ (the class of bounded special functions that will be introduced below).
More precisely we study the deviation bounds for 
$$P_x\left (| \int_0^tf(X_s)ds |\geq t^{\frac12 + \eta} x \right)  $$
for $\eta > 0 $ in the positive recurrent case, and 
in the general null-recurrent case, including the non-regular case, of
$$ P_x \left( | \int_0^t f(X_s) ds | \geq v^{\frac12 + \eta}_t x \right) .$$
Thus we are interested in deviation inequalities in the moderate deviations regime. 
We will
work in the most general situation,
our results hold for any recurrent Markov process taking values in any Polish state space, without restrictions on the 
quality of recurrence. As a counterpart of this general approach, we have to restrict attention to the
class of test functions which are the bounded centered specials functions. In the case of
strong Feller processes, that means that we consider centered functions being of compact support.  

Our approach is based on the regeneration method. Regeneration times
allow to split the trajectory of the process into i.i.d. excursions.
In the one-dimensional case, such regeneration times are usually
introduced as successive visits to recurrent points. Our aim is,
however, to work in the frame of a general state space. In this
general frame, recurrent points usually do not exist, and we use
Nummelin splitting in continuous time, as developed in L\"ocherbach
and Loukianova \cite{dashaeva}, in order to overcome this
difficulty. This technique allows to introduce a recurrent atom for the
process on an enlarged state space and thus to mimic the idea of recurrent points in higher dimensions. Even if the technical details concerning Nummelin
splitting are somehow cumbersome, the evident advantage of this
method is that it is conceptually very easy (regeneration means that
we are able to work with i.i.d. variables) and that it works in any
dimension and for any state space. 

For positive recurrent processes $X,$ we obtain the following deviation upper bound.
For any bounded special function $f$ with $\mu (f) = 0,$ we have   
\begin{eqnarray}\label{eq:mdp1}
P_\pi \left( | \int_0^t f(X_s ) ds | \geq  t^{\frac12 + \eta}   x   \right) &\le &  C_1 \exp \left( - C_2  t^{ 2 \eta}   \;  (x^2 \wedge x)  \right)
+ R_t (x) ,
\end{eqnarray}
where $\pi$ is any initial measure, for any $0 < \eta \le \frac12 .$ Here, $C_1$ and $C_2$ are
explicit positive constants. $R_t(x) $ is a remainder term which is explicitly known and which is of order $ O (  \exp \left( - \sqrt{t} (x \vee 1)   \right)) .$ 
In the general null-recurrent case, including the non-regular case, we have 
\begin{eqnarray}\label{eq:mdp2}
P_\pi \left( | \int_0^t f(X_s ) ds | \geq  v(t)^{\frac12 + \eta}   x   \right) &\le &  C_1 \exp \left( - C_2  v(t)^{  \eta}   \;  (x^2 \wedge x)  \right)
+ R_t (x),
\end{eqnarray}
$R_t$ being of order $ O (  \exp \left( - \sqrt{v(t)} (x \vee 1)   \right)).$ 
Note that (\ref{eq:mdp1}) and (\ref{eq:mdp2}) yield Gaussian concentration in $x.$ This is due to the fact that we consider functions $f$ which are bounded and special. For this class of functions, the contribution of the additive functional up to the first regeneration time can be controlled uniformly in the the starting point of the process, and it is this uniform dependence that is the basis of the Gaussian behaviour. In particular, for uniformly contracting Markov processes, Gaussian concentration is obtained for any bounded function $f$ that need not be of compact support. 

(\ref{eq:mdp1}) and (\ref{eq:mdp2}) represent a first and important step towards the study of moderate deviations in the fully general case
of any recurrent process. The problem of obtaining such non-asymptotic bounds is of major
importance for many applications. Let us cite just some of them :
model selection or other non asymptotic problems for statistics of
Harris processes, particle approximations of Gibbs measures, .... Our work is strongly motivated by applications
to statistics (see for example \cite{driftestimation}), 
in particular model selection, and it is important for such 
applications to obtain bounds that are valid for fixed $t$ such that the
constants involved are as explicit as possible. This is what we try to
achieve in this paper. Note also that our results hold in any dimension and under any starting measure $P_\pi .$

Our deviation inequalities are stated in the moderated deviations regime. Note
that the moderate deviations for additive functionals 
have been widely studied firstly in the case of discrete time, i.e.~for Markov chains. 
We refer the reader
to the work of Guillin (\cite{arnaud0}), Djellout and Guillin (\cite{arnaud1}), Chen and Guillin (\cite{arnaud2})
for a survey on the subject. In the case of time continuous observations, less results are known.
In the ergodic situation, Guillin and Liptser (\cite{arnaud3}) have studied the moderate deviations of 
$ \int_0^t f(X_s) ds $
where $X$ is a multidimensional ergodic diffusion and $f$ any centered function belonging to $L^1 (\mu ) .$  
They use techniques of the stochastic calculus which are well-adapted to the particular diffusion case.
Douc, Fort and Guillin (\cite{DFG}) seem to be the first who deal with the general case of bounded additive 
functionals in the continuous time situation, but under assumptions implying sub-exponential ergodicity of the process (and in particular, positive recurrence) and assuming irreducibility of some skeleton chain (i.e.~assuming minorization of $P_m $ for some $m$ on a petite set). All these results have been achieved in the positive recurrent case, specifying the
quality of ergodicity -- whereas in the present paper we work in the fully general recurrent but not necessarily
positive recurrent case.

Finally, let us cite a recent work by Loukianova et al.~\cite{dasha}
dealing with the same type of questions in the framework of positive recurrent diffusions in
dimension one. In this case regeneration times are hitting times,
and a precise control of these hitting times is provided in \cite{dasha}.
These kinds of results are not available within the general frame
of the present article.

The paper is organized as follows. Since we are using heavily the
method of Nummelin splitting and the concept of deterministic
equivalent, section 2 gives a review of all known results concerning
this technique that will be needed in the sequel. Section 3 gives
the main results, for both positive and null-recurrent cases (Theorem \ref{theo:alltogether}). 
In the case of
null-recurrent but regular models, i.e. in the case when $v_t \sim
t^\alpha$ for some $0 < \alpha < 1,$ some finer control of the
Laplace transform of the length of the regeneration time yields
better results than in the general null-recurrent case (see Theorem \ref{th:dritteversion}).
Sections 4 and 5 are devoted to the proofs. Here, we
also state a technical proposition which is a sort of generalization of
Kac's formula, but more cumbersome in the context of Nummelin
splitting (Proposition \ref{prop:momentp}). Finally,
section \ref{section:galves} gives an idea of possible applications of our results in the
framework
of some interacting particle systems.

\section{Notation}
Consider a probability space $(\Omega, {\cal A}, (P_x)_x) .$ Let $X
= (X_t)_{t \geq 0 }$ be a process defined on $(\Omega , {\cal A},
(P_x)_x) $ which is strong Markov, taking values in a locally
compact Polish space $(E,{\cal E}), $ with c\`adl\`ag paths. $X_0 =
x $ $P_x -$almost surely. We write $L$ for the generator and
$(P_t)_t $ for the semi group of $X$ and we suppose that $X$ is
recurrent in the sense of Harris, with invariant measure $\mu ,$
unique up to multiplication with a constant. This means that for any
set $A \in {\cal E}$ such that $\mu (A) > 0, $ $ \lim\sup_{t \to
\infty} 1_A (X_t ) = 1 $ almost surely. Moreover, we shall write
$({\cal F}_t)_t$ for the filtration generated by the process.

We impose the following condition on the transition semi-group $(P_t)_t $ of $X: $

\begin{ass}\label{regularitypt}
There exists a sigma-finite positive measure $\Lambda $ on $(E,
{\cal E}) $ such that for every $t > 0,$ $P_t (x,dy ) = p_t(x,y)
\Lambda (dy) ,$ where $(t,x,y) \mapsto p_t (x,y) $ is jointly
measurable.
\end{ass}

\begin{ex}
\begin{enumerate}
\item
In general, it is difficult to check whether a given Markov process
is recurrent or not, and the most used criterion for recurrence is
the existence of a so-called Lyapunov-function for the generator of
the process, see for example Meyn and Tweedie (\cite{mt1},
\cite{mt2}, \cite{mt3}).

We say that $V \in dom (L) $ is a Lyapunov-function, if $V \geq 1$
and if there exists a constant $a > 0 ,$ a constant $b$ and a closed
petite set $C$ such that for all $x,$
$$ L V (x) \le - a V(x) + b 1_C (x) .$$

The existence of Lyapunov-functions implies exponential ergodicity. This
concept can be extended to obtain slower rates of convergence, still
in the positive recurrent case, see Douc, Fort and Guillin \cite{DFG}.

\item
In the context of interacting particle systems, (see Galves et al. (\cite{galves}
and section \ref{section:galves}  below), 
recurrence can be shown via arguments using the dual process of the system. 
\end{enumerate}
\end{ex}

\subsection{On the deterministic equivalent of additive functionals}
In the general recurrent not necessarily positive recurrent case, rates of convergence of additive functionals are
given by what is called {\it deterministic equivalent} of additive functionals. This object has been introduced 
 for Markov chains
by Chen \cite{chen} and then been generalized to the context of 
continuous time diffusion models by Loukianova and Loukianov \cite{dashaoleg}
and to any continuous time recurrent Markov process by L\"ocherbach and 
Loukianova in \cite{dashaeva}. In the context of one dimensional diffusion models, similar ideas to
the notion of {\it deterministic equivalent} have also been 
developed in Delattre et al. (\cite{delattre}).
We start by resuming the most relevant results of (\cite{dashaeva})
on the deterministic equivalent that will be needed in the sequel. We first recall the notion of additive
functionals. 

\begin{defin}
An additive functional of the process $X$ is a $\bar{\RR}_+ -$valued, adapted process $A = (A_t)_{t \geq 0} $ such that
\begin{enumerate}
\item
Almost surely, the process is non-decreasing, right-continuous, having $A_0 = 0. $
\item
For any $s, t \geq 0, $ $A_{s+ t } = A_t + A_s \circ \theta_t $ almost surely. Here, $\theta $ denotes the shift operator.
\end{enumerate}
\end{defin}

Examples for additive functionals are $A_t = \int_0^t f(X_s) ds $ where $f$ is a positive
measurable function. Such an additive functional is said to be {\it integrable}, if $ \mu (f) < \infty .$ The deterministic equivalent of any integrable additive functional is a deterministic function $v \mapsto v(t)$ such that $v(0) = 0, v(.) $ is non-decreasing and $v(t) \to \infty $ as $t \to \infty .$ It satisfies that for any integrable additive functional $A_t,$ $A_t / v(t) $ is bounded and bounded away from zero in probability. 
In order to define the deterministic equivalent, we have to recall the notion of a special function (see also \cite{revuz}, \cite{mihai}):
\begin{defin}\label{specialdef}
A measurable function $f : E \to \RR_+ $ is called special if for all bounded and
positive measurable functions $h$ such that $\mu (h) > 0, $ the function
$$ x \mapsto E_x \int_0^\infty \exp \left[ - \int_0^t h(X_s ) ds \right] f(X_t) dt $$
is bounded.
\end{defin}
Note that in the case of strong Feller processes having locally
compact Polish state space, any bounded function having compact
support is special.

By \cite{dashaeva}, any special function $g$ of $X$ with $\mu (g) > 0 $ defines a version of the deterministic equivalent via
\begin{equation}\label{deteq1}
v(t) = E_\pi \int_0^t g(X_s) ds ,
\end{equation}
for any arbitrary initial measure $\pi .$ $v(t)$ is called deterministic equivalent due to the following result (corollary 2.19 of \cite{dashaeva}).

\begin{theo}
For any additive functional $A$ of the process having $E_\mu (A_1) \in ] 0, \infty [, $ for any initial measure
$\pi,$ we have
$$ \lim_{ M \to \infty } \lim \inf_{t \to \infty } P_\pi \left( \frac{1}{M } \le \frac{1}{v(t)} A_t \le M \right) = 1 .$$
\end{theo}

\begin{rem}
The deterministic equivalent is unique up to a constant : for two choices of the deterministic equivalent, $v$ and $v'$, we have that $ \lim_{t \to \infty }  v(t)/v'(t) = c ,$ where $c$ is a positive constant.
\end{rem}

In the sequel, depending on the situation, we shall fix a suitable
choice of $v(t) .$ In the positive recurrent case, evidently $v(t) =
t , $ up to multiplication with a constant. In order to avoid too cumbersome notation, we sometimes also
write $v_t = v(t) .$

\subsection{On Nummelin splitting in continuous time}
The proof of the deviation inequality is based on a very simple idea
: the use of regeneration times that allow to divide the trajectory
of the process into (almost) i.i.d.~excursions. In the
one-dimensional case, regeneration times are introduced as
successive visits to recurrent points. For Harris recurrent Markov
processes with general state space, points are in general not
recurrent. That is why we have to use the Nummelin splitting in
continuous time, as developed in \cite{dashaeva}, in order to
introduce a {\it recurrent atom} for the process. An {\it atom} is a
set that, roughly speaking, behaves as a point for the process. Once
a recurrent atom exists, we can introduce regeneration times that
split the trajectory of the process into (almost) i.i.d. excursions.

We recall briefly the construction of Nummelin splitting in
continuous time.

Introduce a sequence $ (\sigma_n)_{n \geq 1} $ of i.i.d.
$exp(1)$-waiting times, independent of the process $X$ itself. Let $T_0 := 0,$ $T_n := \sigma_1 + \ldots
+ \sigma_n$ and $\bar{X}_n := X_{T_n} .$  Then the chain $\bar{X}= (\bar{X}_n)_n$ is recurrent in the sense of Harris
and its one-step transition kernel $U^1 (x, dy) := \int_0^{\infty} e^{-t} P_t (x, dy) dt $
satisfies the minorization condition
\begin{equation}\label{minoration}
 U^1 (x, dy) \geq \alpha 1_C (x) \nu (dy) ,
\end{equation}
where $0 < \alpha < 1, $ $\mu (C) > 0 $ and $\nu $ a probability measure equivalent to $\mu
(\cdot \cap C)$ (cf \cite{revuz}, \cite{ams}, proposition 6.7). The set $C$ can be chosen to be compact.
Note that the above minorization holds always under the only condition of Harris recurrence. We do not
impose any further condition, neither irreducibility of some skeleton chain nor aperiodicity of the process. 

\begin{rem}\label{choiceofc}
In some cases, the measure of assumption 2.1 $ \Lambda $ satisfies $\Lambda  \sim \mu $, and the densities $p_t (x,y) $ are explicitly known, for example in the case of $k-$dimensional 
diffusions, under suitable regularity assumptions. If one can specify some set $C$ and some time interval $[s, t] $ such 
that 
$$ \inf_{(x,y) \in C \times C, u \in [s, t]} p_u (x,y) > 0, \; \Lambda (C) > 0,$$
(w.l.o.g. also $ \Lambda (C) \le 1 $) then (\ref{minoration}) holds true with 
$$ \alpha = \left[  e^{-s} - e^{-t} \right] \Lambda (C)  \left( \inf_{(x,y) \in C \times C, u \in [s, t]} p_u (x,y) \wedge 1 \right)  \mbox{ and }
\nu = \Lambda ( \cdot \cap C) / \Lambda (C) .$$
In particular, for multi-dimensional diffusions satisfying H\"ormander's condition on some set $\Gamma ,$ the classical results of Kusuoka and Stroock \cite{kuso} allow us to conclude that any choice of a compact set $C \subset \Gamma  $ will be possible.     
\end{rem}

Then it is possible to define on an extension of the original space $(\Omega, {\cal A}, (P_x))$ a Markov process $Z = (Z_t)_{t \geq 0},$ taking values in $E \times [0, 1 ] \times E$ such that the $T_n$ are jump times of the process and such that under $P_x,$ $((Z_t^1)_t, (T_n)_n)$ has the same distribution as $((X_t)_t, (T_n)_n).$ Here are the details of this construction.

First of all, define the following transition kernel $Q ((x,u), dy) $ from $ E \times [0, 1 ] $ to $E :$
\begin{equation}\label{Q}
Q((x,u), dy) = \left\{
\begin{array}{ll}
\nu(dy) & \mbox{ if } (x,u) \in C \times [0, \alpha]\\
\frac{1}{1 - \alpha} \left( U^1 (x, dy) - \alpha \nu(dy) \right)  & \mbox{ if } (x,u) \in C \times ] \alpha , 1] \\
U^1 (x,dy) & \mbox{ if } x \notin C
\end{array} \right. .
\end{equation}

We now recall the construction of $Z_t = (Z_t^1, Z_t^2, Z_t^3) $ taking values in $E\times [0, 1] \times E  $ as given in \cite{dashaeva}. Write $u^1 (x, x') := \int_0^\infty e^{-t} p_t (x, x') dt .$ Let $Z^1_0  =  X_0 = x  .$ Choose $Z_0^2 $ according to the uniform distribution $U$ on $[0, 1] .$ On $\{ Z_0^2 = u\},$ choose $Z_0^3 \sim Q((x,u), dx' ) .$ Then inductively in $n \geq 0,$ on $Z_{T_n} = (x,u,x') :$
\begin{enumerate}
\item
Choose a new jump time $\sigma_{n+1} $ according to
$$e^{-t} \; \frac{p_t(x,x')}{u^1 (x, x') } \, dt  \mbox{ on } \RR_+,$$
where we define $0/0 := a/ \infty :=  1,$ for any $a \geq 0,$
and put $T_{n+1} := T_n + \sigma_{n+1}.$
\item
On $\{ \sigma_{n+1} = t \},$ put $Z_{T_n +s}^2 := u ,$ $Z_{T_n +s}^3 := x'  $ for all $0 \le s < t .$
\item
For every $s < t,$ choose
$$Z_{T_n + s}^1 \sim \frac{p_s(x,y) p_{t-s}(y , x')}{p_t (x,x') } \;  \Lambda (dy ) .$$
Choose $ Z_{T_n + s}^1 := x_0$ for some fixed point $x_0 \in E$ on $\{ p_t (x, x') = 0 \}.$
Moreover, given $Z_{T_n + s}^1 = y,$ on $s + u < t, $ choose
$$ Z^1_{T_n + s+ u } \sim \frac{p_u (y, y') p_{t-s-u}(y' , x')}{p_{t-s} (y, x') } \Lambda (dy') .$$
Again, on $\{ p_{t-s} (y, x')  = 0 \},$ choose $ Z^1_{T_n + s+ u } = x_0 .$
\item
At the jump time $T_{n+1},$ choose $Z^1_{T_{n+1}}  := Z^3_{T_n} = x' .  $ Choose $Z_{T_{n+1}}^2 $ independently of $Z_s, s < T_{n+1}  ,$ according to the uniform law $U.$ Finally, on $\{  Z^2_{T_{n+1}} = u' \},$ choose $ Z^3_{T_{n+1}} \sim Q((x',u'), dx'' ) .$
\end{enumerate}

Note that by construction, given the initial value of $Z$ at time $T_n,$ the evolution of the process $Z^1$  during $[T_n, T_{n+1}[$ does not depend on the chosen value of $Z^2_{T_n} .$

We will write $P_{\pi}$ for the measure related to $X$, under which $X$ starts from the initial measure $\pi(dx)$, and $\P_{\pi}$ for the measure related to $Z$, under which $Z$ starts from the initial measure $\pi(dx)\otimes U(du)\otimes Q((x,u),dy)$. In the same spirit we denote
$E_{\pi}$ the expectation with respect to $P_{\pi}$  and $\E_{\pi}$ the expectation with respect to $\P_{\pi}$. Moreover, we shall write $\FF $ for  the filtration generated by $Z,$ $\GG $ for the filtration generated by the first two coordinates $Z^1 $ and $Z^2 $ of the process, and $\FF^X $ for the sub-filtration generated by $X $ interpreted as first coordinate of $Z.$

Write
$$ A := C \times [0, \alpha ] \times E. $$
$A$ is the recurrent atom of the process. Now we put
$$ S_0 := 0, \; R_0 := 0 , S_{n+1} := \inf \{ T_m > R_n:  Z_{T_m} \in A  \}, R_{n+1} := \inf\{ T_m : T_m > S_{n+1} \}  .$$
Then the sequence of $\FF -$stopping times $R_n$ generalizes the notion of life-cycle
decomposition in the following sense.

\begin{prop}\label{iid} [Proposition 2.13 of \cite{dashaeva}]\\
a) $Z_{R_n + \cdot}$ is independent of ${\cal F}_{S_{n}-} $ for all $ n \geq 1.$\\
b) $Z_{R_n}   \sim \nu(dx) U(du) Q((x,u), dx') $ for all $n \geq 1 .$\\
c) The sequence of $(Z_{R_n})_{n\geq 1} $ is i.i.d.
\end{prop}

\begin{prop}\label{wonderful} [Proposition 2.20 of \cite{dashaeva}]\\
Let $A_t$ be any integrable additive functional of $X.$ Then, up to multiplication by a constant, for any initial measure $\pi$ and any $ n \geq 1,$
$$\E_\pi (A_{R_{n+1}} - A_{R_n}) =  \E_{\nu} (A_{R_1}) = E_\mu (A_1 ) .$$
\end{prop}

From now on, we shall fix a version $\mu $ of the invariant measure such that always
\begin{equation}\label{eq:fixedmu}
 \mu (f) = \E_\pi  \int_{R_1}^{R_2} f(X_s) ds .
\end{equation}

Moreover, we have the following:

\begin{prop}\label{martingale} [Proposition 4.4 of \cite{dashaeva2}]\\
Let $f$ be a measurable $\mu -$integrable function. Put 
$$ \xi_n := \int_{R_{n-1}}^{R_n} f(X_s)ds  , n \geq 1.$$
Then the sequence $(\xi_n)_n$ is a stationary ergodic sequence under $\P_\nu. $ Moreover, for $n \geq 2, $ $\xi_n $ is independent of
${\cal F}_{R_{n-2}}.$
\end{prop}

\begin{rem}
In the usual one-dimensional case, regeneration times $R_n$ allow to split in such a way that
$\xi_n$ is independent of $\xi_{n - 1} .$ Our situation, however, is more complicated, since by construction,
at a regeneration time, $Z_{R_n} $ depends on the state of the process one jump time before, i.e. on $Z_{S_n} .$
This is due to the structure of continuous time and due to the use of Markov bridges (step 3. of the construction).
\end{rem}

Finally, let us recall the following useful result.

\begin{prop}\label{prop:cvonf} [Proposition 2.16 of \cite{dashaeva}]\\
Let $f$ be a special function of the process $X.$ Then
$$ C(f) := \sup_{x \in E} \E_x \int_0^{R_1} |f| (X_s) ds < + \infty .$$
\end{prop}

\section{The deviation inequality}
In this section, we state our main result which is a deviation inequality for
$ \int_0^t f(X_s) ds ,$ where $f$ is a special function of the process having $\mu (f) = 0 .$ Write
$$ N_t := \sup \{ n : R_n \le t \} .$$
In the following, we shall also use the following version of
the deterministic equivalent. 
\begin{equation}\label{eq:fixedvt}
v^*_t = \E_\nu (N_t) +1 ,
\end{equation}
where $\nu = {\cal L} (Z^1_{R_n} ) $ is given in (\ref{minoration}). The
following is our main result.

\begin{theo}\label{theo:alltogether}
Let $f$ be a bounded special function such that $\mu (f) = 0 .$ 
Recall that $ C(f) = \sup_x \E_x \int_0^{R_1} |f|(X_s) ds $ 
and put
$$K(f) := ||f||_{\infty} + C(f), \; B (f) := max( K^2 (f) , K(f) ).$$
1. Suppose that $X$ is positive recurrent. Let 
$$ \hat F (\lambda ) = \E ( e^{- \lambda (R_2 - R_1 ) }) $$
be the Laplace transform of the length of a life cycle and put
$ m = \E( R_2 - R_1)  .$
Then for any $u < m ,$
\begin{equation}\label{eq:legendre}
 \Lambda^* ( u) = \sup_{ \lambda > 0 } \left[- \lambda u  - \log \hat F (\lambda )   \right] > 0,\end{equation}
 and we have the following results. 
\begin{itemize}
\item[(i)]
For any $ 0 < \eta \le \frac12, $ for all $x,$ for any initial measure $\pi  $ and for all $t > 4 m ,$
\begin{eqnarray}
&&P_\pi \left( | \int_0^t f(X_s ) ds | \geq  \, t^{\frac12 + \eta}  \left(\frac{2}{m}\right)^{\frac12 + \eta} x  \right) \nonumber \\
& &\quad \quad \quad \quad \le   4 \exp \left( - t^{ 2 \eta}  \frac{1}{42 m \; B (f)  } \;  (x^2 \wedge x)   \right)
+ 4 e  \exp \left( - \frac{2^{\frac12 + \eta}}{6 \,K (f)  m^{\frac12 + \eta} }t^{\frac12 + \eta}    x  \right) \nonumber \\
&&\quad \quad \quad \quad \quad  +  8   \exp \left( { - t  ( x \vee 1) \frac{3}{4m}  \; \Lambda^* ( \frac23 m ) } \right) .
\end{eqnarray}
\item[(ii)]
For all $x$ and for any initial measure $\pi  ,$ for all $t > 4 m ,$
\begin{eqnarray}
&&P_\pi \left( | \int_0^t f(X_s ) ds | \geq   \, \sqrt t \,  x \frac{\sqrt2}{\sqrt m}  \right) \le   4 \exp \left( -  \frac{1}{42 \; B (f) } \;  (x^2 \wedge x)  \right)
\nonumber \\
&& \quad \quad \quad + 4 e  \exp \left( - \frac{\sqrt 2}{6 \,K(f)   \sqrt m } \sqrt t    x  \right) +  8   \exp \left( { - t ( x \vee 1) \frac{3}{4m} \; \Lambda^* ( \frac23 m ) } \right)  .
\end{eqnarray}
\end{itemize}
2. Suppose that $X$ is null-recurrent. Then there exists $t_0,$ such that for all $ t \geq t_0,$ for all $x, $ for any $ 0 < \eta \le \frac12 ,$ for any initial measure $\pi  $ and for $v^*_t $ as in (\ref{eq:fixedvt}),
\begin{eqnarray}
&&  P_\pi \left( | \int_0^t f(X_s ) ds | \geq (v^*_t)^{\frac12 + \eta}  x  \right) \le
 4 \exp \left( - \frac{1}{42 \; B (f) } \; (v^*_t)^{ \eta} (x^2 \wedge x)  \right)  \nonumber \\
&& \quad \quad \quad \quad \quad + 4  e  \exp \left( - \frac{1}{6  K (f) }(v^*_t)^{\frac12 + \eta}   x  \right) +  8  \exp \left( {-\frac12 (v^* _t)^{\eta } (x\vee 1)  } \right)   .
\end{eqnarray}
$t_0$ is defined through $v^*_{t_0} = 1 .$ 
\end{theo}

The bounds obtained in the above theorem are in general not optimal. But they are valid in a non-asymptotic framework. 
The quantities which are involved are almost all known or known up to some constants. Note also that as a function of the deviation level $x,$ the above bounds $\exp ( - C t^{2 \eta} (x^ 2 \wedge x) )$ and $ \exp ( - C (v_t^*)^{ \eta } (x^ 2 \wedge x ))$ are {\it of Gaussian type}. This is due to the fact that we consider special functions allowing for exponential moments even in cases where the underlying process is {\it not exponentially ergodic}, see also Proposition \ref{prop:momentp} below.

\begin{rem}
In (\ref{eq:fixedvt}) we do not use the usual form of the deterministic equivalent as 
chosen in (\ref{deteq1}). But we have the following comparison result : 
\begin{enumerate}
\item
For the choice of $v_t = E_\nu \int_0^t g(X_s) ds $ as in (\ref{deteq1}), we have 
\begin{equation}\label{eq:upperbound}
 v_t \le C (g) + \mu (g)  \cdot v^*_t  .
 \end{equation}
This can be seen as follows. Since $g > 0,$ we have clearly that
\begin{eqnarray*}
&v_t &\le \sum_{n \geq 0} \E_\nu \left( 1_{ \{ R_n \le t \}} \int_{R_n}^{R_{n+1}}  g(X_s) ds \right) \\
& & \le C(g)  + \sum_{n \geq 1} \E_\nu \left( 1_{ \{ R_{n-1} \le t \}} \int_{R_n}^{R_{n+1}}  g(X_s) ds \right) .
\end{eqnarray*} 
Using Markov's property with respect to ${\cal F}_{R_{n-1}} $ in the last expression and (\ref{eq:fixedmu}), we obtain that 
$$ v_t \le C(g) + \mu (g) [\E_\nu (N_t) +1] ,$$
and thus,  
$$ v_t \le C(g) + \mu (g) \cdot  v_t^* .$$ 
\item
Moreover, we have the following lower bound.
\begin{equation}\label{eq:lowerbound} 
v_t \geq  \mu (g) \cdot v_t^*  - 2 C(g) . \end{equation}
Indeed,
\begin{eqnarray*}
&v_t &\geq \left( \sum_{n \geq 0} \E_\nu \left( 1_{\{ R_n \le t \}} \int_{R_{n+1}}^{R_{n+2}} g(X_s) \right) \right) 
- \E_\nu (\int_{R_{N_t}}^{R_{N_t+2}} g(X_s) )  \\
&& \geq [\E_\nu (N_t) +1]  \mu (g) - 2 C(g) = v_t^* \cdot \mu (g) - 2 C(g),
\end{eqnarray*} 
since $ \E_\nu (\int_{R_{N_t}}^{R_{N_t+2}} g(X_s) ) \le 2 C(g) .$ 
\end{enumerate}
\end{rem}

The above Theorem \ref{theo:alltogether} can be improved in the
null-recurrent but regular case. Here, in accordance with Theorem 3.15 of H\"opfner and L\"ocherbach, \cite{ams},
we call a process {\it regular}, if for $0 < \alpha < 1 $ and a function $l$ 
varying slowly at $\infty, $ the following is true. For any measurable and positive function $g$
with $ 0 < \mu (g) < \infty,$ we have regular variation of resolvants 
\begin{equation}\label{eq:regres}
R_{1/t} g (x ) = E_x \left( \int_0^\infty e^{ - \frac1t s} g(X_s) ds \right) \sim t^\alpha \frac{1}{l(t)} \mu (g) \mbox{ as } t \to \infty ,
\end{equation}
for $\mu-$almost all $x.$ Here we do not consider the case $\alpha = 1 .$ 
We will show in proposition \ref{prop:reg} below that (\ref{eq:regres}) is equivalent to the following :

\begin{equation}\label{eq:domainofattraction}
\P( R_2 - R_1 > x) \sim \frac{x^{- \alpha} l(x)}{\Gamma (1 - \alpha)
} \mbox{ as } x \to \infty .
\end{equation}

Thus we are in the situation where $R_2 -
R_1 $ belongs to the domain of attraction of a stable law. In this case
 a finer control of $N_t$ is possible and it is well-known that
$$ v(t) \sim t^\alpha
\frac{1}{l(t)} \mbox{ as } t \to \infty ,$$
see for instance H\"opfner and L\"ocherbach \cite{ams}, theorem 5.6.A. Moreover,
\begin{equation}\label{eq:fixedchoicevt3}
 1 - \hat F (\lambda)  \sim \lambda^\alpha l( 1/  \lambda ) 
  \mbox{ as } \lambda \to 0  ,
\end{equation}
see for instance Bingham et al. (\cite{bgt}, corollary 8.1.7).

Then we get the following version of the deviation theorem.

\begin{theo}\label{th:dritteversion}
Suppose that (\ref{eq:domainofattraction}) holds.
Let $f$ be a bounded special function such that $\mu (f) = 0 .$ Then there exists $t_0 \geq 0 $ which is given explicitly in (\ref{eq:t}) below,  a function $L$ varying slowly at infinity, such that for all $x $ and for any initial measure $\pi  ,$
for all $ t \geq t_0, $ for all $ 0 < \eta \le \frac{\alpha }{2} ,$ 
\begin{eqnarray*}
&&P_\pi \left( | \int_0^t f(X_s ) ds | \geq t^{\frac{\alpha}{2} + \eta }  x  \right)
\le 4   \exp \left( - \frac{1}{42 \; B (f) } \; t^{ 2 \eta  / (2 - \alpha)}(x^2 \wedge x)  L(t)
  \right)  \\
 &&\quad \quad \quad \quad \quad \quad \quad \quad \quad + 
 4  e  \exp \left( - \frac{1}{6 K (f)  } t^{ \frac{\alpha}{2} + \eta}     x  \right)  +  8  \exp \left( - \frac12  t^{ 2 \eta  / (2 - \alpha)} (x \vee 1)  \, \Lambda^*_t \right)  ,
\end{eqnarray*}
where 
$$ \Lambda^*_t = \sup_{\lambda > 0 } \left[ - \log \hat F (\lambda ) t^{\alpha - \alpha \frac{2 \eta}{2 - \alpha }} \frac{1}{L(t)} - \lambda t^{ 1 - \frac{2 \eta }{2 - \alpha}} \right] > 0 $$
is positive and does not depend on $t $ asymptotically :
$$ \lim_{t \to \infty} \inf \Lambda_t^* \geq (1 - \alpha ) \alpha^{ \alpha/ (1 - \alpha)} > 0 .$$ 
\end{theo}

\begin{rem}
Note that in the above theorem, for $ \eta = \alpha /2 ,$ we obtain a rate of decay of the order of $\exp( - t^{ \alpha / (2 - \alpha)}).$ For $\alpha = 2 ,$ this gives the rate
of convergence $\exp (-  t^{1/3 } ) $ which is better than the rate $\exp (- t^{1/4})$ obtained in Theorem \ref{theo:alltogether}, item 2. (for $\eta = \frac12$). 
\end{rem}

The above Theorems \ref{theo:alltogether} and \ref{th:dritteversion} do not use the Legendre transform of $ \xi_n = \int_{R_{n-1}}^{R_n} f(X_s) ds $ since in general 
the law of the abstract regeneration times $R_n $ is not explicitly known. That is why the above deviation inequalities involve the constant $K(f) .$
In the following, we show how to compare $K(f)$ and $\mu (f) .$

\begin{rem}
\hspace{0,01cm}
\begin{enumerate}
\item
By proposition 2.16 of \cite{dashaeva} we know the following. Recall the definition of the set $C$ of (\ref{minoration}), see also
remark \ref{choiceofc}.
Let $S$ be the first jump time of a Poisson process having rate $1_C (X_s) ,$ and let $K$ be a constant such that for all $x , $
$$ E_x \int_0^S |f| (X_s) ds  \le K , \quad |f(x) | \le K  . $$
Then we have that 
$$ C(f) \le K + \frac{3 K}{\alpha } \mbox{ and hence } K(f) \le 2 K+ \frac{3 K}{\alpha }  . $$
\item
We are now going to explore the relationship between $C(f)$ and the invariant measure $\mu (|f |) $ in some special cases. 

Suppose that the process $X$ is strong Feller and positive recurrent. 
Let $F$ be a compact set such that the support of $f$ is contained in $F.$ Let 
$$ T_C := \inf \{ t \geq 0 : X_t \in C \} $$
be the entrance time in the set $C.$  Recall that the measure $\nu$ is concentrated on $C.$ 
Then 
\begin{equation}\label{eq:cf}
C(f) = \sup_{x \in F} \E_x \int_0^{R_1} |f| (X_s) ds \le  \left[ \| f\|_\infty \sup_{x \in F} E_x T_C \right] + 
\sup_{x \in F \cap C} \E_x  \int_0^{R_1} |f| (X_s) ds .
\end{equation}
By continuity of the map $x \mapsto  \E_x \int_0^{R_1} |f| (X_s) ds ,$ there exists $x_0 \in C \cap F $ such that
$$ \sup_{x \in F \cap C} \E_x  \int_0^{R_1} |f| (X_s) ds  =\E_{x_0}  \int_0^{R_1} |f| (X_s) ds .$$
Moreover, there exists $\varepsilon > 0 ,$ such that for all $x \in B_\varepsilon (x_0), $ 
$$ \E_x \int_0^{R_1} |f| (X_s) ds  \geq \frac12 \E_{x_0}  \int_0^{R_1} |f| (X_s) ds .$$ 
Then we have that 
\begin{eqnarray}\label{eq:cf2}
& \mu (|f|) = \E_\nu \int_0^{R_1} |f| (X_s) ds &\geq \int_{B_\varepsilon (x_0)} \nu (dx ) \E_x \int_0^{R_1} |f| (X_s) ds \nonumber
\\
&&\geq \frac12 \cdot \nu (B_\varepsilon (x_0)) \cdot \E_{x_0}  \int_0^{R_1} |f| (X_s) ds .
\end{eqnarray}
Putting together (\ref{eq:cf}) and (\ref{eq:cf2}), we conclude that 
\begin{equation}\label{eq:cf3} 
C(f) \le  \left[ \| f\|_\infty \sup_{x \in F} E_x T_C \right] + \frac{2}{\nu (B_\varepsilon (x_0))} \mu (|f|) .\end{equation}
 \item
Suppose that $X$ is a positive recurrent one-dimensional diffusion process
$$ d X_t = b(X_t) dt + \sigma (X_t) d W_t .$$
In the case of dimension one, we can avoid Nummelin splitting since 
successive visits of recurrent points allow to split the trajectory into i.i.d. excursions. More
precisely, let $a < b  $ be two recurrent points of $X$ and define a sequence of stopping times $(S_n)_n,
(R_n)_n$ as follows. $S_0 = R_0 = 0,$ 
\begin{equation}\label{eq:regenerationtimes}
 S_1 = \inf \{ t \geq 0 : X_t = b \} , \; R_1 = \inf \{ t \geq S_1 : X_t = a \} , 
 \end{equation}
and for any $n \geq 1 , $ $S_{n+1} = R_n + S_1 \circ \theta_{R_n}, $ 
$ R_{n+1} = R_n + R_1 \circ \theta_{R_n} .$
Using this sequence $(R_n)_n,$ under additional regularity assumptions (see \cite{driftestimation} for the details), we have that 
$$ C(f) \le \kappa \mu (|f|) ,$$
for any function $f$ having compact support. Here, the constant $\kappa $ is explicitly known. 
\end{enumerate}
\end{rem}

\section{Proof of Theorem \ref{theo:alltogether}}
The following proposition will be crucial in the sequel. It states that additive functionals built from bounded special functions
admit a certain number of exponential moments. 

\begin{prop}\label{prop:momentp}
Let $f$ be a bounded special function. Recall that $K(f) = || f||_{\infty} + C(f).$ Then we have for any initial measure $\pi$ and any $n \geq 1, $ that
\begin{equation}\label{eq:momentp2}
 \E_\pi  (|\xi_n|^p) \le p! K(f)^p .
 \end{equation}
In particular, we obtain for any $0 < \lambda < K(f)^{-1}, $
$$   \E_\pi (e^{\lambda \xi_n}) \le 1 + \sum_{p \geq 1} \lambda^p K(f)^p = \frac{1}{1 - \lambda K(f)},$$
and if $\mu (f) = 0,$ 
$$   \E_\pi (e^{\lambda \xi_n}) \le 1 + \sum_{p \geq 2} \lambda^p K(f)^p = 1 + \frac{\lambda^2 K(f)^2 }{1 - \lambda K(f)} . $$
\end{prop}

\begin{proof}
Evidently, we have that
\begin{eqnarray}\label{eq:momentp3}
\xi_n^p &=& \int_{R_{n-1}}^{R_n}\ldots \int_{R_{n-1}}^{R_n} f(X_{t_1}) \ldots f(X_{t_p}) d t_1 \ldots d t_p \\
&\le & p! \int_{R_{n-1}}^{R_n}\ldots \int_{R_{n-1}}^{R_n} 1_{\{ t_1 \le \ldots \le t_p\}}| f|(X_{t_1}) \ldots |f|(X_{t_p}) d t_1 \ldots d t_p.\nonumber
\end{eqnarray}
Taking expectation and conditional expectation with respect to ${\cal F}_{t_{p- 1}} ,$ we get
\begin{eqnarray}\label{eq:successive}
&& \E_\pi \int_{R_{n-1}}^{R_n} \ldots  \int_{R_{n-1}}^{R_n} 1_{\{ t_1 \le \ldots \le t_p\}}| f|(X_{t_1}) \ldots |f|(X_{t_p}) d t_1 \ldots d t_p \nonumber \\
&&=  \E_\pi  \int_{R_{n-1}}^{R_n} \ldots  \int_{R_{n-1}}^{R_n} 1_{\{ t_1 \le \ldots \le t_{p-1}\}}| f|(X_{t_1}) \ldots |f|(X_{t_{p-1}}) \nonumber \\
 && \quad \quad \quad \quad \quad \quad \quad \quad \quad \E_{ Z_{t_{p-1}}} \left[ \int_0^{R_1} |f| (X_s) ds \right] d t_1 \ldots d t_{p-1} .
\end{eqnarray}
But note that for any fixed $z ,$ for $M = ||f||_\infty ,$
\begin{eqnarray*}
\E_z \left[ \int_0^{R_1} |f| (X_s) ds \right]& \le & M \E_z (T_1) + \E_z \int_{T_1}^{R_1} |f| (X_s) ds \\
& \le & M \E_z (T_1) + \E_z ( \E_{Z_{T_1}^1} \int_0^{R_1} |f| (X_s) ds ) \\
& \le & M \E_z (T_1) + C(f) .
\end{eqnarray*}
(Compare to proposition 2.16 of \cite{dashaeva}.) But by construction, for $z = (x,u,x') ,$
$$ \E_z (T_1) = \int_0^\infty t e^{-t} \frac{p_t (x,x')}{u^1 (x, x') } dt ,$$
and this expression does only depend on $x$ and $x' .$ According to proposition 4.1 of \cite{dashaeva2}, we have that
$$ {\cal L} (Z^3_t | Z_t^1 = x) (dx') = u^1 (x, x') \Lambda (dx') .$$
Taking now conditional expectation in (\ref{eq:successive}) with respect to ${\cal F}_{t_{p- 2}} $ and putting all these results together, we obtain
\begin{eqnarray*}
&& \E_\pi \int_{R_{n-1}}^{R_n} \ldots  \int_{R_{n-1}}^{R_n} 1_{\{ t_1 \le \ldots \le t_p\}}| f|(X_{t_1}) \ldots |f|(X_{t_p}) d t_1 \ldots d t_p \nonumber \\
&&=  \E_\pi \int_{R_{n-1}}^{R_n} \ldots  \int_{R_{n-1}}^{R_n} 1_{\{ t_1 \le \ldots \le t_{p-1}\}}| f|(X_{t_1}) \ldots |f|(X_{t_{p-1}})
\E_{ Z_{t_{p-1}}} \left[ \int_0^{R_1} |f| (X_s) ds \right] d t_1 \ldots d t_{p-1}\\
& & \le  \E_\pi \int_{R_{n-1}}^{R_n} \ldots  \int_{R_{n-1}}^{R_n} 1_{\{ t_1 \le \ldots \le t_{p-1}\}}
| f|(X_{t_1}) \ldots |f|(X_{t_{p-1}}) \left[M  \E_{ Z_{t_{p-1}}}  (T_1) + C(f) \right] d t_1 \ldots d t_{p-1}\\
&& = \E_\pi \int_{R_{n-1}}^{R_n} \ldots  \int_{R_{n-1}}^{R_n} 1_{\{ t_1 \le \ldots \le t_{p-2}\}}| f|(X_{t_1}) \ldots |f|(X_{t_{p-2 }}) \\
&& \quad \quad \quad \quad \quad \quad \quad \quad \quad \E_{ Z_{t_{p-2}}} \left[ \int_0^{R_1} |f| (X_s) ( M \E_{Z_s} (T_1) + C(f))  ds  \right]d t_1 \ldots d t_{p-2 }\\
&& = \E_\pi  \int_{R_{n-1}}^{R_n} \ldots  \int_{R_{n-1}}^{R_n} 1_{\{ t_1 \le \ldots \le t_{p-2}\}}| f|(X_{t_1}) \ldots |f|(X_{t_{p-2 }}) \\
&& \quad \quad \quad \quad \quad \quad \quad \quad \quad \E_{
Z_{t_{p-2}}} \left[ \int_0^{R_1} |f| (X_s) ( M  + C(f))  ds
\right]d t_1 \ldots d t_{p-2 },
\end{eqnarray*}
where the last equality follows from
\begin{eqnarray*}
\E_z \int_0^{R_1}  |f| (X_s)  \E_{Z_s} (T_1) ds & = & \E_z \int_0^\infty 1_{\{ s < R_1 \}} |f| (X_s) \left(\int u^1( X_s, x') \Lambda (dx') \E_{(X_s, x')} (T_1)  \right) ds  \\
& = &  \E_z \int_0^\infty 1_{\{ s < R_1 \}} |f| (X_s) ds ,
\end{eqnarray*}
since
$$ \int u^1 (x, x') \Lambda (dx') \int_0^\infty t e^{-t} \frac{p_t (x,x')}{u^1 (x, x') } dt = 1 .$$
Taking successively conditional expectations with respect to ${\cal F}_{t_{p- 3}},  \ldots , {\cal F}_{t_1}$ yields the result.
\end{proof}

A first step in order to prove Theorem \ref{theo:alltogether} is the following proposition. Let 
$$ v_t = \left\{ 
\begin{array}{ll}
\frac{2}{m} t& \mbox{ if $X$ is positive recurrent}\\
v_t^* & \mbox{ if $X$ is null-recurrent}
\end{array}
\right\} . $$

\begin{prop}\label{bernstein}
Let $f$ be a bounded special function such that $\mu (f) = 0 .$ 
Then there exists $t_0,$ such that for all $x > 0 , $ for any initial
measure $\pi  ,$ for any $ 0 \le \eta \le \frac12,$ $0 \le  \delta \le  2 \eta ,$ for any fixed
choice of a deterministic equivalent $v_t  $ and for all $t \geq t_0 ,$ 
\begin{eqnarray*}
 P_\pi \left( | \int_0^t f(X_s ) ds | \geq v_t^{\frac12 + \eta}  x \right)& \le&   4 \exp \left( - \frac{1}{42} \;\frac{x^2 \wedge x}{max( K^2 (f) , K(f))}\; 
 v_t^{ 2 \eta - \delta}  \right) \\
 &&+ 4 e \;  \exp \left( -
\frac{1}{6 K (f)   }v_t^{\frac12 + \eta }   x   \right) +  4 \P_\pi ( N_t > v_t^{1 + \delta} (x \vee 1) )  .
\end{eqnarray*}
Here, $t_0$ is given by the equation $ v_{t_0}= 1 . $
\end{prop}

\begin{proof}
Write $\xi_n := \int_{R_{n-1}}^{R_n} f(X_s) ds ,$ then the $\xi_n, n
\geq 2,$ are identically distributed random variables having mean
zero, such that $ \xi_n $ and $\xi_{n +2} $ are independent. It can
be proven that for any $p \geq 1, $ $\E (|\xi_n|^p) < \infty ,$ by
the properties of a special function. More precisely, we have that
for any $n \geq 1, $ 
\begin{equation}\label{eq:momentp}
 \E (|\xi_n|^p) \le p! K(f)^p,
 \end{equation}
which will be shown in proposition \ref{prop:momentp} below.
Thus for $\lambda > 0,$ sufficiently small ($ \lambda < K(f)^{-1}$ suffices),
\begin{equation}\label{eq:finiteexpmoment}
 Z( \lambda ) := \E_\pi ( e^{\lambda \xi_n } ) , \; n \geq 2,
\end{equation}
exists and is finite (and does not depend on $n$). 

Since the $\xi_n$ are not independent, but only $2-$independent, we
have to proceed in the following way. Firstly, define a sequence
$\xi^{(1)}_n $ by
\begin{equation}\label{eq:xi1}
 \xi^{(1)}_n = \left\{
\begin{array}{ll}
\xi_n & \mbox{ if $n$ odd} \\
0 & \mbox{ elseif }
\end{array}
\right\} .
\end{equation}
Then define a second sequence $\xi_n^{(2)} $ by
$$ \xi^{(2)}_n = \left\{
\begin{array}{ll}
\xi_n & \mbox{ if $n$ even} \\
0 & \mbox{ elseif }
\end{array}
\right\} .$$

Note that $\xi_n$ is not independent of ${\cal F}_{R_{n-1}} ,$ but it is independent of 
${\cal F}_{R_{n-2}} .$ That is why we introduce the following two sub-filtrations, 
associated to the sum of odd and the sum of even terms. 
Let
$$ {\cal G}_n^{(1)} := \sigma \{ R_k, \xi_k^{(1)} : k \le n , k \mbox{ odd } \} , $$
and 
$$ {\cal G}_n^{(2)} := \sigma \{ R_k, \xi_k^{(2)} : k \le n , k \mbox{ even } \} .  $$
Moreover, let 
$$ N_t^{(1)} := \sup \{ n : n \mbox{ odd } , R_n \le t \}, \; N_t^{(2)} := \sup \{ n : n \mbox{ even } , R_n \le t \} .$$
Then it is immediate that $ N_t^{(1)} +2 $ is a $ ( {\cal G}^{(1)}_n)_n -$ stopping time and 
$ N_t^{(2)} +2 $ a $ ( {\cal G}^{(2)}_n)_n -$ stopping time. Moreover, we evidently have that 
$N_t^{(1)} \le N_t, $ $ N_t^{(2) } \le N_t .$ 

Now, for $\lambda $ sufficiently small, let
$$M^1_n := \exp ( \lambda \sum_{k = 2}^n \xi^{(1)}_k  ) \cdot  Z(\lambda)^{-[(n-1) /2] } , \;
M^2_n := \exp ( \lambda \sum_{k = 2}^n \xi^{(2)}_k  ) \cdot
Z(\lambda)^{-[n /2] } ,$$ where $[.] $ denotes the integer part of a
real number. Then $(M^1_n)_n$ and $(M^2_n)_n$ are discrete
$ {\cal G}^{(1)}_{n}-$martingalges  ($ {\cal G}^{(2)}_{n}-$martingales, respectively). Hence
using Doob's stopping rule for positive super-martingales we get
\begin{equation}\label{eq:doob}
\E_\pi \left[  \exp( \lambda \sum_{k=2}^{N^{(1)}_t +2 } \xi_k^{(1)} ) \exp ( - [(N^{(1)}_t +1) /2] \log Z( \lambda ) ) \right] \le 1 
\end{equation}
and also
\begin{equation}
 \E_\pi \left[  \exp( \lambda \sum_{k=2}^{N^{(2)}_t +2 } \xi_k^{(2)}  ) \exp ( - [(N^{(2)}_t+2) /2] \log Z( \lambda ) ) \right] \le 1 . 
\end{equation}
Now, we proceed as follows. Evidently,
$$  \P_\pi(  \int_0^{R_{N_t +2 }} f(X_s) ds \geq v_t^{\frac12 + \eta}  x  ) \le I_0 + I_1 +I_2 ,$$
where $I_0 =\P_\pi ( \xi_1 \geq v_t^{\frac12 + \eta } \, x  /3 ) $,
$$ I_1 = \P_\pi ( \sum_{ k=2}^{N^{(1)}_t +2 } \xi_k^{(1)} \geq v_t^{\frac12 + \eta} x /3 ),\;
 I_2 = \P_\pi ( \sum_{ k=2}^{N^{(2)}_t +2 } \xi_k^{(2)} \geq v_t^{\frac12 + \eta} x  /3 ) .
$$
We start with a study of the first term $I_0 .$ Let
$$ \tilde Z (\lambda ) = \E_\pi e^{\lambda \xi_1 } .$$
Then, for $0 < \lambda < K(f)^{-1} ,$ 
\begin{eqnarray*}
&&  \P_\pi [ \xi_1 \geq  v_t^{\frac12 + \eta }   x  /3 ]   \le \exp \left( - (v_t^{\frac12 + \eta}   \lambda  x/3
  - \log \tilde Z( \lambda ))    \right) .
\end{eqnarray*}
But due to (\ref{eq:momentp}), we have that $ \log \tilde Z(\lambda ) \le \frac{\lambda  K (f)}{1 - 
\lambda K(f) } ,$ and thus, taking $\lambda = \frac12  K(f)
^{-1},$
\begin{equation}\label{eq:i0}
I_0 =  \P_\pi [ \xi_1 \geq v_t^{\frac12 + \eta } x K(f) /3 ] \le e \; \exp \left( -
\frac{1}{6 K (f)   }v_t^{\frac12 + \eta }   x   \right)  .
\end{equation}
Let us now turn to the study of $I_1$ and $I_2 .$ 

Note that for any $k > 0,$
$$ I_1 \le  \P_\pi ( \sum_{ k=2}^{N^{(1)}_t +2 } \xi_{k}^{(1)} \geq v_t^{\frac12 + \eta }  x  /3 ; \; N_t \le k ) + \P_\pi (N_t > k ) .$$
By (\ref{eq:momentp}), we have that for any $0 < \lambda < K(f)^{-1}, $
$$ Z(\lambda )  \le 1 + \sum_{n \geq 2} [ \lambda K(f) ]^n . $$
Hence
$$ Z(\lambda ) \le 1 + \frac{\lambda^2 K(f)^2}{1 - \lambda C(f)}, \mbox{ thus  }
 \log Z( \lambda ) \le   \frac{\lambda^2 K(f)^2}{1 - \lambda K(f) }  .$$

Hence for any $0 < \lambda < K(f)^{-1}, $  using (\ref{eq:doob}) and recalling that $N_t^{(1)} \le N_t , $
\begin{eqnarray*}
&&  \P_\pi [  \sum_{ k=2}^{N^{(1)}_t +2 } \xi_{k}^{(1)} \geq v_t^{\frac12 + \eta }  x 
/3 \; ; N_t \le k ] \le  \P_\pi [  \sum_{ k=2}^{N^{(1)}_t +2 } \xi_{k}^{(1)} \geq v_t^{\frac12 + \eta }  x 
/3 \; ; N^{(1)} _t \le k ]
  \\
  &&\le  \P_\pi \left[  M^1_{N^{(1)}_t +2}  \geq \exp( \lambda v_t^{\frac12 + \eta } x /3  - \log Z(\lambda ) [(N^{(1)}_t+1)/2] ) \; ; N^{(1)} _t \le k \right] \\
 &&\le  \P_\pi \left[  M^1_{N^{(1)}_t +2}  \geq \exp( \lambda v_t^{\frac12 + \eta } x /3  - \frac{\lambda^2 K(f)^2}{1 - \lambda K(f) } [(N^{(1)}_t+1)/2] ) \; ; N^{(1)} _t \le k\right] \\
&&\le  \P_\pi \left[  M^1_{N^{(1)}_t +2}  \geq \exp( \lambda v_t^{\frac12 + \eta } x /3  - \frac{\lambda^2 K(f)^2}{1 - \lambda K(f) }[(k+1)/2] ) \right] \\
&& \le \exp \left( - (v_t^{\frac12 + \eta }  \lambda x \frac{1}{3} - \frac{\lambda^2 K(f)^2}{1 - \lambda K(f) } k)  \right) .
\end{eqnarray*}
In the same way we get that
$$ \P_\pi [  \sum_{ k=2}^{N^{(2)}_t +2 } \xi_k^{(2)} \geq v_t^{\frac12 + \eta }  x  /3  \; ; N_t \le k ]  \le \exp \left( - (v_t^{\frac12 + \eta }  \lambda x \frac{1}{3} - 
\frac{\lambda^2 K(f)^2}{1 - \lambda K(f) } k)  \right)
.$$
Now, take
\begin{equation}\label{eq:choiceofk}
k = v_t^{1 + \delta } \, (x \vee 1)  .
\end{equation}
Then we have that
\begin{equation}\label{eq:zitieren}
 \P_\pi [  \sum_{ k=2}^{N^{(1)}_t +2 } \xi_k^{(1)} \geq v_t^{\frac12 + \eta }  x  /3  \; ; N_t \le v_t^{1 + \delta } ]  \le \exp \left(-v_t^{1 + \delta}   h(t, x) ) \right),
\end{equation}
$$ \P_\pi [  \sum_{ k=2}^{N^{(2)}_t +2 } \xi_k^{(2)} \geq v_t^{\frac12 + \eta }  x  /3  \; ; N_t \le v_t^{1 + \delta } ]  \le \exp \left(-v_t^{1 + \delta}  h(t, x) ) \right),$$
where
\begin{equation}\label{eq:h}
 h(t, x) := (x \vee 1) \sup_{0 <  \lambda < K(f)^{-1} } \left( \frac{\lambda (x \wedge 1)}{3 v_t^{\frac12 + \delta - \eta }} -  \frac{\lambda^2 K(f)^2}{1 - \lambda K(f) }  \right).
\end{equation} 
It can be shown, see for example Birg\'e and Massart (\cite{birge}),
lemma 8, pages 366 and 367, that
\begin{equation}\label{bm}
 \sup_{0 <\lambda < 1/v} \left( \lambda y - \frac{\lambda^2 v^2}{1 - \lambda v}\right) \geq \frac{y^2}{2 v y + 4 v^2} .
 \end{equation}
This is seen as follows. 
A simple calculus shows that 
$$\sup_{0 <\lambda < 1/v} \left( \lambda y - \frac{\lambda^2 v^2}{1 - \lambda v}\right) = \lambda^* y - \frac{(\lambda^*)^2 v^2}{1 - \lambda^* v} ,$$
where 
$$ \lambda^* = \frac{1}{v} \left( 1 - \sqrt{\frac{v}{y + v}} \; \right) < \frac1v  .$$ 
It follows that 
$$\lambda^* y - \frac{(\lambda^*)^2 v^2}{1 - \lambda^* v} = \left( \sqrt{\frac{1}{y + v}} - 1 \right)^2 = \frac{y^2}{ yv + 2 v^2 + 2 v^2 (1 + \frac{yv}{v^2})^{1/2} } \;  ,$$
and using $ (1 + y ) ^{1/2} \le 1 + y/2,$ one gets the desired inequality.

Using (\ref{eq:h}) and (\ref{bm}), we get for any $t \geq t_0$ such that $ (v_{t_0})^{\frac12 + \delta - \eta } \geq 1 ,$ 
\begin{equation}\label{eq:h2}
 h(t, x) \geq  \frac13 \,(x \vee 1) \,   \frac{(x^2 \wedge 1) v_t^{- 1 + 2 \eta -2 \delta} }{2 K(f)v_{t}^{-(\frac12 + \delta - \eta) }+ 12 K^2 (f) }  \geq \frac{1}{42} \;  \frac{( x^2 \wedge x) v_t^{- 1 + 2 \eta -2 \delta}}{ max(K^2(f) , K(f))} .
\end{equation}
Note that the right hand side of (\ref{eq:h2}) tends to zero at speed $v_t^{- 1 + 2 \eta -2 \delta} .$ 
However, this right hand side has to be multiplied with $ v_t^{ 1 + \delta }, $ compare to (\ref{eq:zitieren}),  which yields the term 
$v_t^{2 \eta - \delta }$ which does not tend to zero since by assumption, $ \delta \le 2 \eta .$ Thus, together with (\ref{eq:i0}), \begin{eqnarray*}
 \P_\pi [
\int_0^{R_{N_t +2 }} f(X_s) ds \geq v_t  x   ] &\le &  2 \exp
\left( - \frac{1}{42} \;  \frac{x^2 \wedge x  }{ max(K^2(f) , K(f))} \; v_t^{2 \eta  - \delta}  \right)  \\
&& + e
\; \exp   \left( -
\frac{1}{6 K (f)   }v_t^{\frac12 + \eta }   x   \right)+  2 \P_\pi (N_t > v_t^{1 + \delta}(x\vee 1)  )  .
\end{eqnarray*}

Applying the same argument to $-f$ instead of $f$ yields
\begin{eqnarray}\label{eq:b1}
 \P_\pi(  |\int_0^{R_{N_t +2 }} f(X_s) ds | \geq v_t  x  )& \le&4 \exp
\left( - \frac{1}{42} \;  \frac{x^2 \wedge x}{ max(K^2(f) , K(f))} \; v_t^{2 \eta  - \delta}  \right)  \\
&& + 2 e
\;  \exp \left( -
\frac{1}{6 K (f)   }v_t^{\frac12 + \eta }   x   \right)+  4 \P_\pi (N_t > v_t^{1 + \delta} (x\vee 1)  )  . \nonumber 
\end{eqnarray}

Moreover, note that
\begin{eqnarray*}
\P_\pi \left( | \int_t^{R_{N_t + 2}} f(X_s) ds | \geq v_t^{\frac12 + \eta} x \right)  & \le  & 
\P_\pi \left( | \int_t^{R_{N_t + 1}} f(X_s) ds | \geq \frac12 v_t^{\frac12 + \eta} x \right)\\
&& + 
\P_\pi \left( | \int_{R_{N_t + 1}}^{R_{N_t + 2}} f(X_s) ds | \geq \frac12 v_t^{\frac12 + \eta} x \right)\\
\end{eqnarray*}
and, using Markov's property with respect to ${\cal F}_t$ and to ${\cal F}_{R_{N_t+1}},$ 
\begin{eqnarray}\label{eq:onedim}
\P_\pi \left( | \int_t^{R_{N_t + 1}} f(X_s) ds | \geq \frac12 v_t^{\frac12 + \eta} x \right) &\le 
&  \E_\pi \E_{X_t} e^{ \lambda  \int_{0}^{R_{ 1}}|
f|(X_s) ds } e^{- v_t^{\frac12 + \eta}   \lambda x /2 } .
\end{eqnarray}
Write
$$\tilde{Z}( \lambda) = \E_\pi  \E_{X_t} (\exp  \lambda \int_0^{R_1} |f| (X_s) ds ) .$$
As before we have that
$$ \log \tilde{Z}( \lambda) \le \frac{\lambda  K(f)}{1 -  \lambda K(f)},$$
and thus,  taking $\lambda  = \frac12 K(f) ^{-1} , $
$$ \P_\pi \left( | \int_t^{R_{N_t + 1}} f(X_s) ds | \geq \frac12 v_t^{\frac12 + \eta}  x \right)   \le  e  \; \exp \left( - \frac{1}{4 K (f) } v_t^{\frac12 + \eta} x \right)  .$$
In the same way we get that 
$$ \P_\pi \left( | \int_{R_{N_t + 1}}^{R_{N_t + 2}} f(X_s) ds | \geq \frac12 v_t^{\frac12 + \eta}  x \right)   \le  e  \; \exp \left( - \frac{1}{4 K (f) } v_t^{\frac12 + \eta} x  \right)  .$$
\end{proof}

\begin{rem}
The most delicate point in the above proof is (\ref{eq:h2}) which gives a lower bound converging to zero. It is important to be able  
control the speed of convergence of this expression to zero. If we decided to work with the Legendre transform, then the above proof remains true
by replacing (\ref{eq:h}) by 
$$ h(t,x) = (x \vee 1) \sup_{\lambda > 0 } \left( \lambda \frac{x \wedge 1 }{v_t^{\frac12 + \delta - \eta }} - \log Z( \lambda ) \right) .$$ 
In this form it is evident that $h(t,x)  \to 0 $ as $t \to \infty ,$ and developing $ Z( \lambda ) =  1 + \frac{\sigma^2 }{2} \lambda^2 + o (\lambda^2)  ,$ where $$  \sigma^2 := \E (\int_{R_1}^{R_2} f(X_s) ds )^2 ,$$
would yield the same speed of convergence. However, using this approach we would not be
able to put hands on the exact form of the constants appearing in the remainder term $ o (\lambda^2),$ and that is why we decided to use the upper bound of $\log Z( \lambda )$ as proposed in the proof. 
\end{rem}

In order to make use of Proposition \ref{bernstein}, we have to control 
the number of life cycles before time $t,$ $N_t .$ We first study the positive recurrent case. Recall that $v_t = \frac{2}{m}t $ in this case. 
\begin{prop}\label{prop:33}
Suppose that $X$ is positive recurrent. Then
$$ \P_\pi ( N_t > v_t (x \vee 1)    )  \le  2 \exp \left(  { - \frac12 v_t (x \vee 1) \sup_\lambda h_t (\lambda , x)}  \right),$$
where
$$ h_t( \lambda , x)  := - \left( 1 - \frac{m}{t ( x \vee 1)} \right) \log  \hat F  (\lambda )  - \lambda \frac{m}{2(x \vee 1)}    $$
and where
$$ \hat F (\lambda ) = \E ( e^{- \lambda (R_2 - R_1 ) }) $$
is the Laplace transform of the length of a life cycle. Moreover we have for any $ u < m $ that 
\begin{equation}
 \Lambda^* ( u) = \sup_{ \lambda > 0 } \left[- \lambda u  - \log \hat F (\lambda )   \right] > 0 .\end{equation}
\end{prop}

\begin{proof} {\bf of Proposition \ref{prop:33}}
Writing $k_t := [v_t (x \vee 1) ],$ for any $\lambda > 0 ,$
$$ \P_\pi (N_t > v_t( x \vee 1) ) = \P_\pi (R_{k_t} \le t ) \le  \P_\pi (R_{k_t} - R_1 \le t ) \le \P_\pi ( e^{ - \lambda (R_{k_t}- R_1)  } \geq e^{ - \lambda t } ) .$$
But, using the definition of $\xi_n^{(1)}$ and $\xi_n^{(2)}$ as in (\ref{eq:xi1}), with $f = 1 ,$ and the same technique as above, 
\begin{eqnarray}\label{eq:gut}
 \P_\pi ( e^{ - \lambda (R_{k_t}- R_1)  } \geq e^{ - \lambda t } ) &\le&
\P_\pi (  e^{ - \lambda \sum_{n=2}^{k_t} \xi_n^{(1)} } \geq e^{ - \lambda  t/2 } ) + \P_\pi (  e^{ - \lambda \sum_{n=2}^{k_t} \xi_n^{(2)} } \geq e^{ - \lambda t/2 } )\nonumber  \\
&\le & 
2  \hat F (\lambda)^{(k_t - 1)/2 } e^{\lambda  t/2 } \le 2 \hat F (\lambda)^{\frac{v_t (x\vee 1) - 2}{ 2}  } e^{\lambda  t/2 },
\end{eqnarray}
where
$$ \hat F (\lambda ) = \E ( e^{- \lambda (R_2 - R_1 ) }) $$
is the Laplace transform of the length of a life cycle. So write
$$ h_t( \lambda , x)  := - \left( 1 - \frac{2}{v_t (x \vee 1)} \right) \log  \hat F (\lambda )  - \lambda \frac{t}{v_t (x \vee 1) }  .$$
Then
$$ \P_\pi (N_t > v_t (x \vee 1)) \le 2  e^{ -\frac12  v_t  (x \vee 1)  \sup_\lambda h_t (\lambda , x)} .$$
Let us finally show that 
$$  \Lambda^* ( u) = \sup_{ \lambda > 0 } \left[- \lambda u  - \log \hat F (\lambda )   \right] > 0 $$
for any $ u < m .$  This is evident using the fact that for $\lambda \to 0,$
$  - \log \hat F  (\lambda ) \sim 1 - \hat F (\lambda ) \sim m \lambda .$
This concludes the proof.
\end{proof}

{\bf Proof of Theorem \ref{theo:alltogether}, item 1. }
The proof follows from Proposition \ref{bernstein} taking $\delta = 0 $ and from 
Proposition \ref{prop:33} in the following way :
A straightforward calculus shows that, for all $t \geq 4m , $ since $x \vee 1 \geq 1,$  
$$ \sup_\lambda h_t (\lambda , x ) \geq \frac34 \sup_\lambda \left( - \frac23 m \lambda - \log \hat F  (\lambda )   \right) = \frac34 \; \Lambda^* ( \frac23 m ) .$$
This concludes the proof of item 1. of Theorem \ref{theo:alltogether}. 
\hfill $\bullet$

In the null-recurrent case, recall the choice of
a deterministic equivalent $v_t^*$ of (\ref{eq:fixedvt}).
Then the following proposition allows us to control the deviations for $N_t$ in 
the null-recurrent case.

\begin{prop}\label{prop:momentdent}
$$\E_\pi (N_t^p) \le p! (v^*_t) ^p .$$
\end{prop}

\begin{proof}
$$ N_t = \sum_{n \geq 1 } 1_{\{ R_n \le t \} } ,$$
and thus
$$ N_t^p = p! \sum_{ 1 \le n_1 \le n_2 \le \ldots \le n_p } 1_{\{ R_{n_1} \le t \} }\cdot \ldots \cdot 1_{\{ R_{n_p } \le t \} }  = p! \sum_{ 1 \le n_1 \le n_2 \le \ldots \le n_p }  1_{\{ R_{n_p } \le t \} } . $$
Hence, using successively Markov's property with respect to ${\cal F}_{R_{n_{p-1}}} $ etc, we get
$$ \E_\pi (N_t^p ) \le  p! \sum_{ 1 \le n_1 \le n_2 \le \ldots \le n_{p-1} } \E_\pi ( 1_{\{ R_{n_{p-1} } \le t \} } \E_{Z_{R_{n_{p-1} }}} ( \sum_{n_p \geq 1 } 1_{\{R_{n_p}\le t \}} ) ) . $$
But note that for any $n,$
\begin{eqnarray*}
 \E_{Z_{R_n}} ( \sum_{k \geq 1 } 1_{\{R_{k}\le t \}} ) ) \le 1 + \E_{Z_{R_n}} ( \sum_{k \geq 2 } 1_{\{R_{k}\le t \}} ) ) &\le& 1 +   E_{Z_{R_n}} ( E_{Z_{R_1}} \sum_{k \geq 1 } 1_{\{R_{k}\le t \}} ) )\\
 & = & 1 + \E_\nu \sum_{k \geq 1 } 1_{\{ R_k \le t \}} = v^*(t) ,
 \end{eqnarray*}
since ${\cal L} (Z_{R_1} ) = \nu .$ Induction on $n$ then yields the assertion.
\end{proof}

We get the following corollary.

\begin{cor}
We have that
\begin{equation}\label{eq:n1}
 \P_\pi( N_t > (v^*_t)^{1 + \delta} (x \vee 1) ) \le 2 \exp \left( {-\frac12 (v^*_t)^\delta  } (x \vee 1)  \right) .
\end{equation}
\end{cor}

\begin{proof}
By proposition \ref{prop:momentdent}, 
$$ \E_\pi ( \exp ( \lambda N_t ) ) \le \frac{ 1}{1 -   \lambda v^*_t },$$
for $\lambda $ sufficiently small. Choosing $\lambda = (v^*_t)^{-1}/2 $, we get that
$$ \P_\pi( N_t > (v^*_t)^{1 + \delta} (x \vee 1)  ) \le e^{- \lambda (v^* _t)^{1 + \delta} (x \vee 1) } \E_\pi ( \exp ( \lambda N_t))\le 2  e^{- \frac12 (v^*_t)^{\delta  } (x \vee 1)} .$$
\end{proof}

{\bf Proof of Theorem \ref{theo:alltogether}, item 2.}
Choosing $\delta = \eta , $ item 2. follows immediately from Proposition \ref{bernstein} and from (\ref{eq:n1}) above. 
\hfill $\bullet$

\section{Proof of Theorem \ref{th:dritteversion}}

\begin{proof} {\bf of Theorem \ref{th:dritteversion}}
The proof is a slight modification of the proof of proposition \ref{bernstein}.
We go back to the proof of proposition \ref{bernstein}. We fix the following choice of $k$ (compare to (\ref{eq:choiceofk})) :
\begin{equation}\label{eq:rightspeed}
k = k_t = 2 + (x \vee 1) t^\gamma \frac{1}{l ( 1 / \lambda_t) }  , \mbox{ where } \gamma = \alpha + 2 \eta  \frac{1 - \alpha}{2 - \alpha} .
\end{equation}
Here, $ \lambda_t \to 0 $ will be defined in (\ref{eq:lt}) below and $l$ is the function of (\ref{eq:fixedchoicevt3}).  

As in (\ref{eq:gut}), we have for any $\lambda > 0,$ 
$$ \P_\pi (N_t > k) \le 2 e^{ - \frac12 \left[  - \log  \hat F (\lambda )  (k-2)   - \lambda t \right] } ,  $$
and therefore, due to the choice (\ref{eq:rightspeed}), we have to find
$$ \sup_\lambda \left[ - \log \hat F (\lambda )(x \vee 1)  t^\gamma \frac{1}{l ( 1 / \lambda_t) }  - \lambda t \right] .$$
Writing $L(t)  = l ( 1 / \lambda_t),$ this last expression can be rewritten as 
$$   t^{ \frac{\gamma - \alpha }{1 - \alpha}}\, (x \vee 1)  \sup_\lambda \left[ - \log \hat F (\lambda )  t^{\alpha \frac{1 - \gamma}{1 - \alpha}} \frac{1}{L(t) }  - \lambda \frac{t^{\frac{1- \gamma}{1 - \alpha }}}{x \vee 1} \right] \geq   t^{ \frac{2 \eta }{2 - \alpha}} (x \vee 1) \, \Lambda_t^* ,$$
where 
$$  \Lambda_t^* =  \sup_\lambda \left[ - \log \hat F (\lambda )  t^{\alpha \frac{1 - \gamma}{1 - \alpha}} \frac{1}{L(t) }  - \lambda t^{\frac{1- \gamma}{1 - \alpha }} \right] .$$ 
(Note that $ \frac{1 - \gamma}{1 - \alpha} = 1 - \frac{2 \eta}{2 - \alpha} .$) 
Write for simplicity 
$$ s_t:=  t^{ \frac{1 - \gamma}{1 - \alpha}} = t^{1 - \frac{2 \eta}{2 - \alpha}} .$$
Then,  
$$ \Lambda_t^* = \sup_{ \lambda > 0 } \left[ - \log \hat F (\lambda)\frac{1}{L(t) } s_t^\alpha - \lambda s_t \right] = \sup_{ \lambda > 0 } \left[ - \log \hat F (\lambda)\frac{1}{l ( 1 / \lambda_t) } s_t^\alpha - \lambda s_t \right] ,$$
and we have to show that $ \Lambda_t^*$ is positive. 
Note that since $R_2 - R_1$ does not possess any moments, we are not able to develop $\hat F (\lambda ) $ near $0 $ in the usual
way. But we can use (\ref{eq:fixedchoicevt3}). That's why we 
 take $\lambda $ of the form $\lambda_t \to 0  $ at a speed that will be precised in (\ref{eq:lt}) below. 
Note that $ \log \hat F (\lambda_t ) = \log  \left[ 1 - (1 - \hat F (\lambda_t)) \right] .$ Due to  (\ref{eq:fixedchoicevt3}) we have that 
$$  - \log \hat F (\lambda_t ) \frac{1}{l ( 1 / \lambda_t) }s_t^\alpha   - \lambda_t s_t   \; \sim \lambda_t^\alpha s_t^\alpha  - \lambda_t s_t,$$
and hence
\begin{equation}\label{eq:dasha1}
\Lambda_t^*  \geq  - \log \hat F (\lambda_t ) \frac{1}{l ( 1 / \lambda_t) } s_t^\alpha   - \lambda_t s_t \; \sim \lambda_t^\alpha \;  s_t^\alpha - \lambda_t s_t\mbox{ as } t \to \infty  .\end{equation}

Maximizing $ \lambda^\alpha \; s_t^\alpha - \lambda s_t $ with respect to $\lambda $ suggests
the choice 
\begin{equation}\label{eq:lt}
\lambda_t =   \alpha^{1/(1 - \alpha)} s_t^{- 1 } .
\end{equation}

For this choice we get that
$$ \lambda_t^\alpha \; s_t^\alpha - \lambda_t s_t  =(1 - \alpha ) \alpha^{ \alpha/ (1 - \alpha)} > 0  . $$
This implies that $\Lambda_t^* $ is strictly positive eventually and that 
$$  \lim_{t \to \infty} \inf \Lambda_t^* \geq (1 - \alpha ) \alpha^{ \alpha/ (1 - \alpha)} > 0 .$$

We continue the proof following the lines of the proof of proposition \ref{bernstein}. (\ref{eq:zitieren}) is true with $v_t^{\frac12 + \eta } $ replaced by $ t^{ \frac{\alpha }{2} + \eta } $ and $ v_t^{1 + \delta } $ by
$k_t   .$ Here, 
$$  h(t,x) = \sup_{\lambda > 0 } \left( \lambda   \frac{x}{3 k_t \, t^{ -\frac{ \alpha }{2} - \eta    }} -  \frac{\lambda^2 K(f)^2}{1 - \lambda K(f) }  \right)
= \sup_{\lambda > 0 } \left( \lambda   \frac{x\wedge 1 }{3 \tilde k_t \, t^{ -\frac{ \alpha }{2} - \eta    }} -  \frac{\lambda^2 K(f)^2}{1 - \lambda K(f) }  \right), $$
where $\tilde k_t = \frac{2}{x \vee 1} + t^\gamma / L(t).$ 
Let 
\begin{equation}\label{eq:t}
t_0 \mbox{ such that for all } t \geq t_0, \frac{ t^\gamma}{L(t)}  t^{ -\frac{ \alpha }{2} - \eta   }\geq 1.
\end{equation}
Then for $t \geq t_0,$ with $B(f) = max ( K^2 (f), K(f)),$ 
$$ h(t,x) \geq \frac13 \; (x^2\wedge 1) \; L(t)^{2  }  \left[\frac{1}{2 (x \wedge 1) K(f)  + 12 K^2 (f)  } \right ] t^{ -2 \gamma +  \alpha + 2 \eta } 
\geq \frac{x^2 \wedge 1}{ 42 B (f) } \;  L(t)^{2  } \;  t^{ -2 \gamma +  \alpha + 2 \eta } .$$

\end{proof}

\begin{rem}
Note that in the regular case, we even have convergence in law of $ (\frac{1}{\sqrt{v_n}} \int_0^{nt} f(X_s) ds)_t $ to $\sigma  B \circ W^\alpha ,$ where 
$B$ is a one-dimensional Brownian motion and $W^\alpha $ the Mittag-Leffler process of index $\alpha ,$ i.e. the process inverse of the stable
subordinator of index $\alpha ,$ see for instance Touati \cite{touati90}.  
\end{rem}

There are various examples where the exact form of the Laplace transform is known. Brownian motion in dimension one is the most famous example. 

\begin{ex}
Let $X $ be the one-dimensional standard Brownian motion. In this case, we can define the regeneration times without
Nummelin splitting in the following very simple way. Let
$$ R_n = \inf \{ t > S_n : X_t = 0 \},  \; S_n = \inf \{ t > R_{n-1} : X_t = 1 \} , \; R_0 = 0 .$$
Then
$$ \hat F ( \lambda ) = e^{ - 2 \sqrt{2 \lambda }} ,$$
see Revuz-Yor (\cite{revuzyor}). In this case it is possible to take
$v_t = \sqrt{t} $ (see for example H\"opfner and L\"ocherbach
(\cite{ams}), theorem 5.6.A.). 
\end{ex}

We conclude this section with the following proposition showing that the notion of regularity (\ref{eq:domainofattraction}) is intrinsic
of the process and does not depend on the concrete splitting we are using. 

\begin{prop}\label{prop:reg}
Suppose that for $0 < \alpha < 1 $ and a function $l$ 
varying slowly at $\infty, $ the following is true. For any measurable and positive function $g$
with $ 0 < \mu (g) < \infty,$ we have regular variation of resolvants 
\begin{equation}\label{eq:regres2}
R_{1/t} g (x ) = E_x \left( \int_0^\infty e^{ - \frac1t s} g(X_s) ds \right) \sim t^\alpha \frac{1}{l(t)} \mu (g) \mbox{ as } t \to \infty ,
\end{equation}
for $\mu-$almost all $x.$ 
Then we have
\begin{equation}\label{eq:domainofattraction2}
\P( R_2 - R_1 > x) \sim \frac{x^{- \alpha} l(x)}{\Gamma (1 - \alpha)
} \mbox{ as } x \to \infty .
\end{equation}
\end{prop}

\begin{proof}
By theorem 3.15 of \cite{ams}, we know that (\ref{eq:regres2}) implies weak convergence
$$ \frac{(A_{tn})_{t\geq 0}}{n^\alpha / l(n)} \to E_\mu (A_1) W^\alpha ,$$
for any additive functional $A_t$ of the process having $0 < E_\mu (A_1) < \infty .$
This convergence holds true in $D(\RR_+ , \RR ) ,$ under $P_\pi $ for any initial measure $\pi .$ 
$W^\alpha $ is the Mittag-Leffler process of index $\alpha .$
 
From now on we fix $\nu $ as initial measure. Interpreting $X$ as first coordinate of $Z$ and using Chacon-Ornstein's theorem, 
we have in particular weak convergence of $(N_{tn}) / ({ n^\alpha / l(n)}), $ under $\P_\nu .$  Now, write $\mbox{\bf v} (t) = t^\alpha / l(t) $ and
let $\mbox{\bf a} (n) $ be its asymptotic inverse, i.e. $  \mbox{\bf a} (n) \sim n^{1/\alpha } \tilde l (n) .$ Note that $   \mbox{\bf a} (n) / n \to \infty $ since $\alpha < 1.$   
$N_t$ is the inverse process of the sequence of regeneration times $R_n,$ which implies that 
\begin{equation}\label{eq:cool}
 \frac{R_n}{  \mbox{\bf a} (n)}  \mbox{ converges weakly as } n \to \infty , \mbox{ under } \P_\nu ,
 \end{equation}
 where the limit is non-degenerated. 
But $R_n = R_1 + (R_2- R_1) + \ldots + (R_n - R_{n-1}) .$ We would like to deduce from this that necessarily (\ref{eq:domainofattraction2})
is true. Unfortunately, due to the complex definition of the Nummelin splitting,
the $(R_k- R_{k-1})_k$ are not independent, but they have all the same law. Independence holds true only after ``nearly'' exponential
times, the reason for this being the choice of the new jump times in step 1. of the splitting algorithm.  

That is why we proceed as follows. Define for any $k \geq 0,$ $\tilde T_k = \inf \{ T_n : T_n > R_k \} $ the first
jump of the process after the $k-$th regeneration time. Then it is straightforward to show that under $\P_\nu,$ $ \tilde T_k - R_k $
is exponentially distributed with parameter $1,$ and that the random times $\tilde T_k - R_k , k \geq 0, $ are independent.  
Hence, 
$$ \frac{\sum_{k = 0}^{n-1} (\tilde T_k - R_k )}{n } \to 1, \mbox{ which implies that }  \frac{\sum_{k = 0}^{n - 1} (\tilde T_k - R_k )}{\mbox{\bf a} (n)} \to 0 $$
almost surely. This together with (\ref{eq:cool}) yields weak convergence of 
\begin{equation}\label{eq:cool2}  \frac{\sum_{k = 0}^{n-1} (R_{k+1} - \tilde T_k  )}{\mbox{\bf a} (n)} .\end{equation} 
But $R_{k+1} - \tilde T_k  , k \geq 0 ,$ is an i.i.d. sequence of random variables under $\P_\nu, $ thus (\ref{eq:cool2}) implies that 
necessarily the law of $ R_{k+1} - \tilde T_k$ belongs to the domain of attraction of a stable law, see for instance Feller, \cite{feller}, XIII.6.  
This implies that 
$$ \P_\nu ( R_2 - \tilde T_1 > x ) \sim x^{- \alpha }  L(x) $$
for some function $L$ varying slowly at infinity. Since 
\begin{eqnarray*}
& \P_\nu ( R_2 - \tilde T _1  > x ) \le  \P_\nu ( R_2 - R_1  > x )& \le \P_\nu ( R_2 - \tilde T_1  > x/2 ) + \P_\nu ( \tilde T_1 - R_1 > x/2)\\
&& = \P_\nu ( R_2 - \tilde T_1  > x/2 ) + e^{ -x/2}   ,
\end{eqnarray*}
we deduce that 
 $$ x \mapsto  \P_\nu ( R_2 - R_1  > x )  \mbox{ varies regularly at infinity with index } \alpha .$$
This implies (\ref{eq:domainofattraction2}). 
\end{proof}

\section{Application to some interacting particle systems}\label{section:galves}
In some applications it is interesting to know at which speed empirical means $\frac1t \int_0^t f(X_s) ds $ converge to the -- in general unknown --
invariant measure $\mu (f) . $ This is most often the case in statistical applications. 

Consider for instance the following interacting particle system which has been studied in Galves et al. \cite{galves}. Particles
are on positions (sites) of a finite subset $V \subset \ZZ^d,$ and any particle has either spin $+1$ or $-1 .$ So let $A = \{ -1, +1\} $
and let $E = A^{V } .$ Any element of $E$ shall be called configuration of the system. Configurations will be denoted by
letters $\eta, \xi , \zeta .$ Any element $i \in V $ is called a site. For any site $i$ let $\eta^i $ be the modified configuration
$ \eta^i (i ) = - \eta (i), \eta^i (j) = \eta (j) $ for all $j \neq i .$ We associate to any site $i$ and any configuration $\eta$ a 
spin-flip-rate $c_i (\eta) \geq 0 .$ Roughly speaking, site $i$ will change its spin at rate $c_i (\eta) $ whenever the overall
configuration of particles is $\eta .$ We suppose that 
$$ \sup_{\eta} c_i (\eta) \le M_i $$ for some constant $M_i .$ 
Then the associated interacting particle system is a Markov process $X$ on $E$ having 
generator 
$$ L f (\eta) = \sum_{i \in V} c_i (\eta) [ f( \eta^i) - f( \eta) ] .$$
Let $V_{i} (k) = \{j \in \ZZ^d; 0 \le \|j
- i \| \le k\},$ where $\|j\| = \sum_{u=1}^d |j_u|$ is the usual
$L_1$-norm of
$\ZZ^d$. In Galves et al. \cite{galves} the following criterion for recurrence of the system has been given :

\begin{prop}[Theorem 3 of Galves et al.]
There exists a sequence $\lambda_i (k ) , k \geq -1 ,$ associated to the spin-flip rates $c_i$ such that :
If 
\begin{equation}\label{subcritical}
 \sup_{i \in V} \sum_{k \geq 0} |V_i (k)| \lambda_i (k) \le 1,
\end{equation}
then the process is recurrent. If the above sum is strictly less than $1,$ then 
the process is uniformly exponentially ergodic. Here, $ |V_i (k)| $ is the number of sites belonging to $ V_i (k) .$
\end{prop}

\begin{rem}
The above theorem relies on the construction of a backward dual process $C_s^{(i)} , s\geq 0,$ such that 
for any site $i \in V,$ $C_s^{(i)}$ denotes the set of sites at time $-s $ that have to be known in 
order to determine the spin of site $ i$ at time $0.$ The cardinal of $C_s^{(i)}$ can
be compared to a classical branching process in continuous time, having reproduction mean 
$\sum_{k \geq 0} |V_i (k)| \lambda_i (k) < 1, $ and thus being subcritical. We refer the reader to Galves et al. (\cite{galves})
for the details. 
\end{rem}

The sequence $\lambda_i (k) $ can be constructed explicitly, we refer the reader to Galves et al. (2008) for the details.   
In \cite{galves}, we were mainly concerned with the issue of perfect simulation of the invariant measure $\mu $ of the 
process, under condition (\ref{subcritical}), based on the precise knowledge of the spin-flip rates $c_i .$ Suppose now, that we
are in the converse situation, observing the process over some time interval $[0, t] ,$ without knowledge of the spin-flip rates, and
that we want to deduce information about the associated invariant measure (and for example on the associated
spin-flip rates, a posteriori). This is the classical situation of statistical inference. 

Note that assumption 2.1 is satisfied and that the process is strongly Feller. 

We apply Theorem \ref{bernstein} with $\eta = \frac12 .$ Since the process is uniformly exponentially ergodic, any constant function is 
special. In particular, we have for any bounded function $f,$ since $f - \mu (f) $ is special, that 
\begin{equation}
 P \left(| \frac1t \int_0^t f(X_s) ds - \mu (f) | > x \right) \le c e^{ - C t (x^2 \wedge x) } + R_t (x) ,
\end{equation}
for some constants $c, C .$ This is a first step on the way towards estimating spin-flip rates and associated interaction schemes
for interacting particles.

\section{Final remarks}
Our results obtained in the positive recurrent case (Theorem \ref{theo:alltogether}, item 1.) 
have to be compared to the results obtained by
Loukianova et al. (\cite{dasha}). The results stated in the present paper
are different from the situation studied in Loukianova et al. (2009)
since we consider centered functions that are 
special, i.e. satisfying $\sup_x \E_x \int_0^{R_1} |f| (X_s) ds <
\infty .$ Additive functionals built of special functions admit exponential 
moments independently of the degree of recurrence of the underlying process, 
and this is a very special feature about special functions that does in general not
hold for other functions. 

Namely, if $f$ is any positive bounded function having compact
support, then $f$ is certainly special, but $\bar f := f - \mu (f) $
in general won't (constants are special only if the underlying process is
uniformly ergodic). 

Note however that our method works for any function $f$
such that exponential moments (\ref{eq:finiteexpmoment}) are finite.
Thus if, for $\lambda $ sufficiently small,
\begin{equation}\label{eq:hmm}
 \E_\pi (e^{\lambda (R_2 - R_1 )} ) + \E_\pi ( e^{ \lambda R_1 })  < \infty ,
\end{equation}
then our result holds also for any function $\bar f = f - \mu (f) ,$ where $f$ is positive, bounded and of compact support. It will be the subject of some futur work to give necessary conditions for (\ref{eq:hmm}).


\begin{thebibliography}{99}
\bibitem{an} Athreya, K.B.,  Ney, P.
A new approach to the limit theory of recurrent Markov chains.
Trans. Am. Math. Soc. 245 (1978), 493-501.
\bibitem{bgt} Bingham, N.H.; Goldie, C.M.; Teugels, Jozef L.
{\it Regular variation.} Encyclopedia of Mathematics and its
applications, Vol. 27. Cambridge etc.: Cambridge University Press.
XIX (1987).
\bibitem{birge} Birg\'e, Lucien; Massart, Pascal.
Minimum contrast estimators on sieves: Exponential bounds and rates of convergence. 
Bernoulli 4, No.3, 329-375 (1998). 
\bibitem{mihai} Brancovan, M. Fonctionnelles additives sp\'eciales des processus r\'ecurrents au sens de Harris. Z. Wahrscheinlichkeitstheor. Verw. Geb. 47 (1979), 163--194.
\bibitem{chen} Chen, X. How often does a Harris recurrent Markov chain recur? Ann. Probab. 27
(1999), 1324--1346.
\bibitem{arnaud2} Chen, X., Guillin, A..
The functional moderate deviations for Harris recurrent Markov chains and applications.
Annales de l'Institut Henri Poincare (Prob. Stat.)  40 (2004), 89--124.
\bibitem{delattre} Delattre, S., Hoffmann, M., Kessler, M.. Dynamics adaptive estimation of a
scalar diffusion. Preprint 762 of Laboratoire de Probabilit\'es et
mod\`eles al\'eatoires, october 2002.
\bibitem{arnaud1} Djellout, H., Guillin, A.. Moderate deviations of Markov Chains with atom.
Stoch. Proc. Appl. 95 (2001), 203-217.
\bibitem{DFG} Douc, R., Fort, G., Guillin, A.,
Subgeometric rates of convergence of $f$-ergodic strong Markov processes. 
Stochastic Processes Appl. 119, No. 3, 897-923 (2009). 
\bibitem{feller} Feller, W.: {\it An introduction to probability theory and its applications. } Vol. 2 Wiley: New York 1971. 
\bibitem{galves} Galves, A., Garcia, N.L., L\"ocherbach, E.(2008): Perfect simulation and finitary coding for multicolor systems with interactions of infinite range. Submitted, arxiv.org/abs/0809.3494.
\bibitem{arnaud0} Guillin, A.. Uniform moderate deviations of functional empirical processes of Markov Chains.
Probability and Mathematical Statistics 20 (2000), 237-260.
\bibitem{arnaud3} Guillin, A., Liptser, R. Examples of moderate deviation principle for diffusion processes.
Discrete and Continuous Dynamical Systems (B) 6 (2006), 803-828.
\bibitem{ams} H\"opfner, R., L\"ocherbach, E.: {\it Limit theorems for null recurrent Markov
processes.} Memoirs AMS {\bf 161}, Number 768, 2003.
\bibitem{kuso} Kusuoka, S.; Stroock, D.
Applications of the Malliavin calculus. III. 
J. Fac. Sci., Univ. Tokyo, Sect. I A 34, 391-442 (1987).

\bibitem{dashaeva} L\"ocherbach, E., Loukianova, D.. On Nummelin splitting for continuous time Harris recurrent Markov processes and application to kernel estimation for multi-dimensional diffusions. Stoch. Proc. Appl. 118 (2008), No. 8, 1301-1321.
\bibitem{dashaeva2} L\"ocherbach, E., Loukianova, D.. The law of iterated logarithm for additive functionals and martingale additive functionals of Harris recurrent Markov processes. Stoch. Proc. Appl. 119 (2009), 2312-2335.
\bibitem{dasha} L\"ocherbach, E., Loukianova, D., Loukianov, O.. Deviation bounds in ergodic theorem for positively recurrent one-dimensional diffusions and integrability of hitting times. Manuscript, 2009. arxiv.org/abs/0903.2405.
\bibitem{driftestimation} L\"ocherbach, E., Loukianova, D., Loukianov, O.. Penalized nonparametric drift estimation in a continuous time one-dimensional diffusion process. Manuscript, 2009. To appear in ESAIM : P $\& $ S.  
\bibitem{dashaoleg} Loukianova, D., Loukianov, O.. Deterministic equivalent of additive functionals of recurrent diffusions and drift estimation. Statistical Inference for Stochastics processes 111 (2008), No 2, 107--121.
\bibitem{mt1} Meyn, S.P., Tweedie, R.L.
Generalized resolvents and Harris recurrence of Markov processes.
Cohn, Harry (ed.), {\it Doeblin and modern probability. Proceedings of the Doeblin conference `50 years after Doeblin: development in the theory of Markov chains, Markov processes, and sums of random variables'.} Providence, RI: American Mathematical Society. Contemp. Math. 149 (1993), 227-250.
\bibitem{mt2} Meyn, S.P., Tweedie, R.L.. A survey of Foster-Lyapunov techniques for general state space Markov processes. In Proceedings of the Workshop on Stochastic Stability and Stochastic Stabilization, Metz, France, June 1993.
\bibitem{mt3} Meyn, S.P., Tweedie, R.L.. Stability of Markovian
processes III: Foster-Lyapunov criteria for continuous-time
processes. Adv. Appl. Probab. 25 (1993), 487--548.
\bibitem{nummelin} Nummelin, E. A splitting technique for Harris
recurrent Markov chains. Z. Wahrscheinlichkeitstheorie Verw. Geb. 43
(1978), 309--318.
\bibitem{numm84} Nummelin, E.: {\it General irreducible Markov chains and non-negative operators.} Cambridge University Press, Cambridge, England, 1984.
\bibitem{revuz} Revuz, D.: {\it Markov chains.} Revised edition. Amsterdam: North Holland 1984.
\bibitem{revuzyor} Revuz, D., Yor, M.: {\it Continuous martingales and Brownian motion.} 3rd edition, Grundlehren der Mathematischen Wissenschaften 293. Berlin: Springer 2005.
\bibitem{touati90} Touati, A.. Loi fonctionnelle du logarithme it\'er\'e pour les processus de Markov r\'ecurrents. Ann. Probab. 18, No.1 (1990), 140-159.
\end{thebibliography}
\end{document}